\documentclass[11pt]{article}
\usepackage{graphicx} 

\usepackage[dvipsnames]{xcolor}
\usepackage{eurosym}
\usepackage{tikz}
\usepackage[linesnumbered,ruled,vlined]{algorithm2e}
\usepackage{amsmath}
\usepackage{amsfonts}
\usepackage{graphicx}
\usepackage{multirow}
\usepackage{textcomp}
\usepackage{pgfplots}
\pgfplotsset{width=10cm,compat=1.9}
\usepackage{longtable}
\usepackage{lscape}
\usepackage{array}
\usepackage{varwidth}
\usepackage[a4paper,top=3cm,bottom=3cm,left=2.5cm,right=2.5cm]{geometry}

\usetikzlibrary{positioning,shapes,shadows,arrows}

\title{Dynamic Pickup-and-Delivery for Collaborative Platforms with Time-Dependent Travel and Crowdshipping}

\author{Sara Stoia\footnotemark[1]~, Demetrio Lagan\`{a}\footnotemark[2]~, Jeffrey W. Ohlmann\footnotemark[3]~}

\date{}
\begin{document}
\maketitle


\renewcommand{\thefootnote}{\fnsymbol{footnote}}

\footnotetext[1]{Department of Economics, Statistics, and Finance \textquotedblleft Giovanni Anania\textquotedblright, University of Calabria, Arcavacata di Rende (CS), Italy. E-mail address: sara.stoia@unical.it}
\footnotetext[2]{Department of Mechanical, Energy and Management Engineering (DIMEG), University of Calabria, Italy. E-mail address: demetrio.lagana@unical.it}
\footnotetext[3]{Department of Business Analytics, University of Iowa, Iowa City, United States. E-mail address: jeffrey-ohlmann@uiowa.edu}

\begin{abstract}
We study a pickup-and-delivery problem that arises when customers randomly submit requests over the course of a day from a choice of vendors on a collaborative e-commerce portal. Based on the attributes of a customer request, a dispatcher dynamically schedules the delivery service on either a dedicated vehicle or a crowdshipper, both of whom experience time-dependent travel times. While dedicated vehicles are available throughout the day, the availability of crowdshippers is unknown a priori and they appear randomly for only portions of the day. With an objective of minimizing the sum of routing costs, piece-rate crowdshipper payments, and lateness charges, we model the uncertainty in request arrivals and crowdshipper appearances as a Markov decision process. To determine an action at each decision epoch, we employ a heuristic that partially destroys the existing routes and repairs them guided by a parameterized cost function approximation that accounts for the remaining temporal capacity of delivery vehicles. Through a set of computational experiments, we demonstrate the improvement of our approach over a myopic approach in several key performance metrics. In addition, we conduct computational experiments demonstrate the impact of inserting wait time in the route scheduling and the benefit of explicitly modeling time-dependent travel times. Through our computational testing, we also investigate the effect of demand management mechanisms that facilitate many-to-one request bundles or one-to-many request bundles on reducing the cost to service requests.

\vspace{0.1cm}
\end{abstract}




\section{Introduction} \label{sec:intro}

\noindent For decades, the advance of information and communications technology has facilitated a trending increase in e-commerce. In the early years of e-commerce, brick-and-mortar businesses still possessed the advantage of immediacy as delivery from a website was expensive, slow, and/or inconvenient. The expanding reach of companies such as Amazon.com has allowed it to provide a market place with a wide variety of goods and the ability to deliver many of them within an hour or two in urban areas and within a day or two practically anywhere else. 

In this study, we examine the strategy of local brick-and-mortar businesses competing with a large e-commerce provider such as Amazon.com through the establishment of an e-platform that aggregates these local businesses' product offerings and provides delivery service. A primary motivation for local businesses to join a collaborative delivery platform is to compete with the distribution power of companies like Amazon.com (Mattioli 2021 \cite{wsjShops}). Pairing an e-shopping platform with a high level of customer service (professional shopper consultation, favorable return policy, etc.) provides local brick-and-mortar companies to facilitate ``shop local" initiatives and stave off a ``retail apocalypse." The global COVID-19 pandemic that stretched from 2020 to 2023 has altered consumer behavior (perhaps permanently) and local brick-and-mortar stores may need to serve as both a showroom (for in-person shopping) and a warehouse, with employees doubling as salespeople and as order-pickers.

We model the delivery service associated with an e-shopping platform of a consortium of local businesses as a stochastic dynamic pickup-and-delivery problem. In this setting, customers at different locations place requests with various businesses thorough a single e-platform. Each request consists of a pickup location (the business fulfilling the request), a delivery location (the customer placing the request), a ready time specifying when the request will be available for pick-up, and a soft deadline outlining the customer's timeliness expectations. As requests arrive randomly throughout the day, the dispatcher must decide how to serve these requests with the available delivery capacity. 

We consider delivery capacity in the form of a fleet of dedicated vehicles supplemented by the presence of crowdsourced couriers which we refer to as crowdshippers. Dedicated vehicles have known availability for day-long shifts while crowdshippers' appearance times are unknown to the dispatcher a priori and their availability duration (declared at appearance to the system) is relatively short. Therefore, a key decision-making element of our problem is whether to service a request with one of the dedicated vehicles or with one of the crowdshippers that is currently available (or may become available in the future). 


Further, we consider request dispatching decisions in the context of time-dependent travel for delivery vehicles. As stated by Vidal et al. \cite{vidal2021arc}, ``an inadequate management of time-dependent travel times $\ldots$ represents the greatest current barrier to the practical adoption of vehicle routing software." We account for traffic patterns by incorporating travel times that depend on the origin and destination as well as the time of day.

To represent the uncertainty in request arrivals and crowdshipper appearances, we model our problem setting as a Markov decision process (MDP) with the objective of minimizing the cost of fulfilling customer requests over an operating day. To identify a policy for the intractable MDP, we employ an approximate dynamic programming (ADP) technique called cost function approximation (CFA) within a destroy-and-repair heuristic. In our CFA, we modify the current epoch cost function via a parameter, tuned via offline simulation, accounting for the delivery vehicles' remaining temporal capacities. As CFA does not require online modeling of the future, it scales well in problem size and facilitates real-time decisions within less than a minute. 



Our work contributes to the literature by being the first to consider a dynamic pickup-and-delivery problem with time-dependent travel times and delivery capacity consisting of dedicated vehicles and randomly-appearing crowdshippers. We provide a methodological contribution by devising destroy-and-repair heuristic leveraging a CFA to address the tradeoff between assigning a request to a dedicated vehicle versus a crowdshipper. Our computational experiments demonstrate the effectiveness of our CFA approach versus a myopic strategy. In addition, we computationally exhibit the value of explicitly accounting for time-dependent travel times. Finally, in our computational experiments we examine the impact of two separate demand management mechanisms that influence the arrival of customer requests to the e-platform. 

The paper is organized as follows. In Section~\ref{sec:lit}, we position our work with respect to the literature. Section~\ref{sec:probform} formally presents the mathematical modeling of the problem as a Markov decision process (MDP). We describe our solution approach utilizing a cost function approximation in Section~\ref{sec:solution}. In Sections~\ref{sec:instances} and \ref{sec:computation}, we outline our data sets and present our solution approach's computational performance, respectively. We conclude with a summary in Section~\ref{sec:concl}. 

\section{Literature Review} \label{sec:lit}

\noindent To frame our contribution, we discuss research related to the two primary characteristics of our problem: dynamic pickup-and-delivery and crowdsourcing. While we focus on research considering a pickup-and-delivery problem with time windows with stochastic and dynamic elements (DPDPTW) for a broader overview of the literature, we refer the reader to the survey on dynamic pickup and/or delivery problems by Berbeglia et al. \cite{berbeglia2010dynamic}, the characterization of dynamic vehicle routing problems by Psaraftis et al. \cite{psaraftis2016dynamic}, and the survey on stochastic and/or dynamic vehicle routing problems by Ritzinger et al. \cite{Ritzinger2016}. 



Mitrovi{\'c}-Mini{\'c} et al. \cite{mitrovic2004waiting} and Mitrovi{\'c}-Mini{\'c} et al.\cite{mitrovic2004double} consider a DPDPTW and incorporate waiting strategies to implicitly accommodate unknown future requests. Vonolfen \& Affenzeller \cite{vonolfen2016distribution} propose methodology to automatically tailor the waiting strategy based on problem characteristics. While our problem context is complicated by having a heterogeneous fleet of dedicated vehicles and crowdshippers, we also allow a vehicle to wait and our computational experiments demonstrate the value of waiting in our problem setting. For a DPDPTW with different priority classes of requests, Ghiani et al. \cite{GHIANI2022} devise a policy function approximation for which parameters specifying the fraction of the fleet reserved for each request class are learned by a neural network.


The dial-a-ride problem (DARP) refers to a special case of the pickup-and-delivery problem involving passenger transportation (often the elderly or disabled) and the related operational considerations. Xiang et al. \cite{xiang2008study} study a dynamic DARP with time windows and a variety of stochastic events, including uncertain travel times. Leveraging probabilistic information about inbound requests provided by outbound requests, Schilde et al. \cite{SCHILDE2011} apply several metaheuristic approaches to a dynamic DARP with time windows motivated by hospital patients. Schilde et al. \cite{schilde2014integrating} extend Schilde et al.\cite{SCHILDE2011} to account for time-dependent stochastic travel times. However, none of these works consider crowdshipping which is a critical element of our setting. While Souza et al.\cite{souza2022bi} consider a dynamic DARP with a heterogeneous fleet, this fleet varies by capacity and rental cost, not by appearance time and temporal availability as crowdshippers do in our problem.

Schrotenboer et al. \cite{schrotenboer2021platform} consider a dynamic stochastic pickup-and-delivery problem that shares the same motivation as our work -- a delivery platform that consolidates requests from local stores. Schrotenboer et al. \cite{schrotenboer2021platform} develop a parameterized CFA based on the Dantzig-Wolfe reformulation of a set packing formulation that assigns requests to vehicles. Once requests are assigned to a vehicle, they cannot be reassigned, requiring an urgency parameter in the CFA that determines whether a request is assigned at an epoch or postponed for consideration at a future epoch. In contrast, we allow reassignment of a request to a different vehicle up until the point at which a vehicle is en route to the request's pickup location, allowing flexibility to make routing adjustments in response to newly-arriving requests. Further, we supplement a fleet of dedicated vehicles with crowdshippers who randomly appear to the dispatching system while Schrotenboer et al. \cite{schrotenboer2021platform} only consider delivery capacity provided by dedicated vehicles (bikes in their application). Third, in our problem the travel time between stores and customer locations is time-dependent and region-dependent.

Restaurant meal delivery is another application domain for dynamic pickup and delivery that bears similarity to our work. To account for the uncertain ready times at the pickup locations resulting from the stochastic nature of restaurant meal preparation, Ulmer et al. \cite{ulmer2021restaurant} employ a CFA in which they introduce a time buffer parameter that alters the assignment of requests to vehicles by effectively penalizing actions that result in deliveries too close to the customer deadlines. They pair this CFA with a strategy of postponing the assignment of requests to drivers. Our study has several differences relative to the work in of Ulmer et al. \cite{ulmer2021restaurant}. First, we supplement a fleet of dedicated vehicles with crowdshippers who randomly appear to the dispatching system while Ulmer et al. \cite{ulmer2021restaurant} only consider delivery capacity provided by dedicated vehicles. Second, we seek to minimize the sum of routing cost, delivery lateness charges, and piece-rate payments to crowdshippers while Ulmer et al. \cite{ulmer2021restaurant} minimize just the delivery lateness. Third, in our problem the travel time between stores and customer locations is time-dependent and region-dependent. While we do not consider uncertain ready times for requests, we do note that if there are time-dependent patterns in the ready times of requests (as may be the case with food preparation), these can be incorporated into the time-dependent travel times. Steever et al. \cite{STEEVER2019} also consider a food delivery problem in which the defining characteristic is that customers may order from more than one restaurant. In our computational experiments, we also consider a set of instances in which there also exists this type of many-pickups-to-one-delivery type of requests. 

The same-day delivery problem (SDDP) refers to a variant of the stochastic dynamic pickup-and-delivery problem in which all orders are picked up from a single origin (distribution center). Of particular relevance is the work of \cite{ulmer2018same} that considers a SDDP with a heterogeneous fleet for which they determine request assignment (truck or drone) with a policy function approximation (PFA) based on geographical proximity to the distribution center. However, this PFA is not directly applicable to our heterogeneous fleet of dedicated vehicles and crowdshippers that differ with respect to temporal capacity and must serve requests with varying pickup locations. Chen et al. \cite{chen2022deep} also consider a SDDP with a heterogeneous fleet of trucks and drones and apply deep reinforcement learning to incorporate more detailed state information in the request assignment policy. While deep RL is a promising methodology, Chen et al. \cite{chen2022deep} consider fleet sizes smaller than in our setting and their results suggest that the benefit of their approach declines as fleet size increases. Dayarian \& Savelsbergh \cite{dayarian2020crowdshipping} consider a dynamic SDDP in which in-store customers complement a company's fleet to serve online orders; they demonstrate the benefit of incorporating probabilistic information about future online orders and in-store crowdshippers versus a myopic approach. Our modeling framework for crowdshippers is flexible enough to accommodate the in-store customer as a special case. Silva et al. \cite{silva2023deep} also consider a SDDP with crowdshippers integrate a value-based function approximation using a neural network with an optimization recourse model to reduce the action space for a scenario, where a scenario consists of a set of delivery requests and whether or not a matching in-store customer exists.

The uncertain appearance time of crowdshippers that provide delivery capacity is a key feature of our dynamic pickup and delivery model. For an overview of operations research literature on crowdsourcing, we refer to surveys by Savelsbergh \& Ulmer \cite{Savelsbergh2022} and Alnaggar et al. \cite{Alnaggar2021}. We discuss treatments of crowdsourcing delivery in dynamic models. 

Arslan et al. \cite{arslan2019crowdsourced} consider a DPDPTW relying on ad hoc drivers, who dynamically announce the origin and destination of planned journeys from which they will detour as specified amount, as well as a dedicated fleet as a backup option. Arslan et al. \cite{arslan2019crowdsourced} formulate the assignment of requests to ad hoc drivers and dedicated vehicles as matching problem that they reoptimize upon the arrival of new information. An exact solution of this matching problem is facilitated by an appearance lead time for the ad hoc driver, narrow time windows of availability (ranging from 10 minutes to 30 minutes), and limits on the number of stops an ad hoc driver will make. While our problem setting is similar, we assume no appearance lead time for crowdshippers, wide availability time windows of one to four hours, and no limits on the number of stops. 

For a dynamic delivery problem in which orders originate from a single depot, \cite{ulmer2020workforce} focus on the scheduling of company vehicles in the presence of unscheduled crowdshippers with uncertain appearance times and availability durations. Using a continuous approximation and a value function approximation, \cite{ulmer2020workforce} produce a set of shifts for scheduled drivers that is expected to achieve a target service level for delivery. Behrendt et al.\cite{behrendt2022prescriptive} also focus on the scheduling of couriers while accounting for ad hoc crowdshippers and present a prescriptive machine learning method identify a set of shifts for scheduled couriers that minimizes total courier payments and penalty costs for expired requests. 


While it is not a pickup-and-delivery problem, there are some dynamic implementations of the vehicle routing problem with occasional drivers related our work. Archetti et al. \cite{archetti2021online} consider a vehicle routing problem in which a regular fleet of vehicles is complemented by a known set of occasional drivers with specified time windows of availability are deployed to serve both static requests known a priori and online requests. Archetti et al. \cite{archetti2021online} handle the dynamic requests with an insertion approach and periodic re-optimization via variable neighborhood search. Mancini \& Gansterer \cite{mancini2022bundle} consider an extension of the vehicle routing problem with occasional drivers (VRP-OD) in which occasional drivers are assigned customer bundles through a bidding system. Our work differs because we assume that both customer requests and crowdshippers appear dynamically, while Archetti et al. \cite{archetti2021online} assumes occasional drivers are known a priori to serve dynamically-arriving requests and Mancini \& Gansterer \cite{mancini2022bundle} consider dynamically-appearing crowdshippers to serve customer demand that is known a priori.

\section{Problem Description and Model Formulation}
\label{sec:probform}

\noindent In this section, we formally describe and formulate the stochastic dynamic pickup-and-delivery problem for our setting. Table~\ref{table:notationTable} provides an overview of the notation we employ in our description.  

We assume that an e-platform receives the customer requests randomly throughout the day. A customer request $r$ corresponds to a tuple of information: the pickup location ($o_r$), the delivery location ($d_r$), the time which the request will be ready ($e_r$), and a soft deadline representing the latest time it should delivered ($l_r$). If the vehicle handling customer request $r$ arrives at location $o_r$ before time $e_r$, the vehicle must wait until $e_r$ before executing the request pick up. If the vehicle handling customer request $r$ delivers the request at location $d_r$ after $l_r$, it incurs a lateness charge.

\begin{table}
\begin{center}

\caption{Overview of notation.}
\begin{tabular}{V{8cm}V{10.5cm}}
\hline
$\mathcal{T}$ & time horizon of request arrivals, $t \in \mathcal{T} = \left[ 0, 1, \ldots, T \right]$ \\ 
$\mathcal{R}_k^p$ & set of in-process requests at epoch $k$ \\
$\mathcal{R}_k^o$ & set of outstanding requests at epoch $k$ \\
$\mathcal{R}_k$ & set of active requests at epoch $k$,  $\mathcal{R}_k = \mathcal{R}_k^p \cup  \mathcal{R}_k^o$ \\
$\mathcal{U}_k$ & set of newly-arriving requests to the system at epoch $k$\\
$o_r$ & pickup location of request $r$\\
$d_r$ & delivery location of request $r$\\
$e_r$ & ready time of request $r$\\
$l_r$ & deadline of request $r$\\
$\mathcal{V}_k$ & set of dedicated vehicles available at epoch $k$ \\ 
$\mathcal{G}_k$ & set of crowdshippers available at epoch $k$ \\ 
$\mathcal{M}_k$ & entire pool of vehicles available at epoch $k$, $\mathcal{M}_k = \mathcal{V}_k \cup \mathcal{G}_k$ \\
$n^{a}_{m}$ & starting location of vehicle $m$\\
$n^{b}_{m}$ & ending location of vehicle $m$ \\
$a_m$ & appearance time of vehicle $m$ at location $n^a_m$ \\
$b_m$ & time by which vehicle $m$ must arrive at location $n^b_m$ \\ 
$\theta_m$ & routing plan of vehicle $m$\\
$f_m$ & arrival time of vehicle $m$ at the first location in $\theta_m$\\
$w_m$ & departure time of vehicle $m$ from the first location in $\theta_m$\\
$s_k$ & pre-decision state\\
$s_k^x$ & post-decision state\\
$\mathcal{X}(s_k)$ & set of possible actions in $s_k$, $x \in \mathcal{X}(s_k)$\\
$c(s_k,x)$ & cost incurred by action $x$ when system in state $s_k$ at epoch $k$\\
$\rho$ & per-delivery fee (\$) \\
$\mu_1$ & travel time multiplier (\$  per minute) \\
$\mu_2$ & lateness multiplier (\$  per minute) \\
$\omega_k$ & set of new customer requests arriving at epoch $k$\\
$\gamma_k$ & set of new crowdshippers appearing at epoch $k$ \\
$\lambda$ & cost multiplier for a vehicle's remaining delivery time \\
$\eta$ & wait time multiplier for a vehicle's remaining delivery time \\
$\tau(i,j,t)$ & travel time between locations $i$ and $j$ when departing at time $t$ \\
\hline
\label{table:notationTable}
\end{tabular}
\end{center}
\end{table}


To provide delivery service, we assume there is a fleet of dedicated vehicles, $\mathcal{V} = \left\{1,\ldots,V \right\}$, and a set of crowdshippers, $\mathcal{G}= \left\{1,\ldots,G \right\}$. Each dedicated vehicle $v \in \mathcal{V}$ is available over a known shift starting and ending at the depot, $n^a_v = n^b_v = 0$ for $v \in \mathcal{V}$. For example, a dedicated vehicle $v$ available for the entire day would have a shift from time $a_v = 0$ to time $b_v = \hat{T} > T$, where $\hat{T}$ is sufficiently large enough to ensure delivery of all requests. Each crowdshipper $g \in \mathcal{G}$ is available to service requests during a time window $[ a_g, b_g ]$ that is unknown to the dispatcher a priori. Upon appearing to the platform at time $a_g$ and declaring their location $n^a_g$, a crowdshipper $g$ communicates their availability to service requests subject to requirement that they can arrive at location $n^b_g$ by time $b_g$. The crowdshipper's destination requirement prevents the dispatcher from assigning a request that would result in a crowdshipper making a delivery in a geographically distant location near the end of its declared time window.

The objective is to minimize the cost of fulfilling customer requests over an operating day. To measure the cost of providing service, we minimize the sum of three components: (1) the total travel cost of the dedicated vehicles and crowdshippers (as function of travel time), (2) the total lateness charge (as a function of the amount of time customer requests are delivered after their soft deadlines), and (3) the total per-delivery fees paid to crowdshippers. Note that we assume a fixed per-delivery fee of $\rho$ for all requests served by a crowdshipper and include the travel cost of the crowdshippers in the objective function to represent the variable fees paid by the e-platform to the crowdshippers to account for differences in the travel demands of requests due to their relative pickup and delivery locations.  

To represent the stochastic nature of the arrival of customer requests and the appearance of crowdshippers, we model the problem as a Markov decision process (MDP) over finite, discrete-time horizon. In the remainder of this section, we describe the components of the MDP: the states, the actions and deterministic transition to the post-decision state, and the exogenous information and stochastic transition to the pre-decision state.  

\subsection{States}

\noindent At a decision epoch $k$, the pre-decision state $s_k$ contains all the relevant and available information to make dispatching decisions. Specifically, a pre-decision state $s_k$ includes the current time of day ($t_k$), the status of each active request, the status of each active vehicle, and the set of newly-arriving requests ($\mathcal{U}_k$). We refer to a customer request for which service (pickup and delivery) has not been complete as an \emph{active} request and let $\mathcal{R}_k$ be the set of active requests at epoch $k$. Each active customer request $r$ corresponds to a tuple of information, $(o_r, d_r, e_r, l_r)$, corresponding to geographical ($o_r, d_r$) and temporal ($e_r, l_r)$ information. We partition the set of active requests into two subsets, $\mathcal{R}_k^p$ and $\mathcal{R}_k^o$. We refer to each request $r \in \mathcal{R}_k^p$ as \emph{in-process} because a vehicle has already picked up request $r$ or is currently en route to $o_r$ to pick it up. We refer to each request $r \in \mathcal{R}_k^o$ as \emph{outstanding} because a vehicle has not begun its service of request $r$. So, $\mathcal{R}_k = \mathcal{R}_k^p \cup \mathcal{R}_k^o$ and $\mathcal{R}_k^p \cap \mathcal{R}_k^o = \emptyset$. We differentiate between $\mathcal{R}_k^p$ and $\mathcal{R}_k^o$ because an active request $r \in \mathcal{R}_k^o$ can possibly be reassigned to another vehicle at a future epoch while an active request $r \in \mathcal{R}_k^p$ cannot.


We represent the entire pool of vehicles available at epoch $k$ as \\$\mathcal{M}_k = \mathcal{V}_k \cup \mathcal{G}_k = 
\left\{ m \in \mathcal{V} \cup \mathcal{G}: t_k \in [ a_m, b_m ] \right\}$. We represent the attributes of vehicle $m \in \mathcal{M}_k$ by the tuple $\left( \theta_m, b_m, f_m, w_m \right)$, where $\theta_m = \langle \theta_m(1), \theta_m(2), \ldots \rangle$ is a sequence of locations composing vehicle $m$'s planned route, $b_m$ is the end time of vehicle $m$'s availability, $f_m$ is the arrival time of vehicle $m$ at the first location in $\theta_m$, and $w_m$ is the planned departure time from the first location in $\theta_m$. If $f_m > t_k$, vehicle $m$ is en route to location $\theta_m(1)$. If $f_m = t_k$, then vehicle $m$ has just arrived at location $\theta_m(1)$ at epoch $k$. If $f_m < t_k$, then vehicle $m$ is idle at location $\theta_m(1)$ with a planned departure time of $w_m$. 

We formally describe a pre-decision state as $s_k = \left( t_k, \mathcal{R}_k^p, \mathcal{R}_k^o, \mathcal{U}_k, \mathcal{M}_k \right)$.  Algorithm \ref{alg:PDS} outlines the dynamics of the Markov decision process. Lines~\ref{line:init}--\ref{line:endwait} describe the deterministic transition from state $s_k$ to post-decision state $s_k^x$ as the result of an action $x$ identified by the CFA-based destroy-and-repair heuristic (DRACE). Lines~\ref{line:incrtime}--\ref{line:stochend} describe the stochastic transition from a post-decision state $s_k^x$ to a pre-decision state $s_{k+1}$ as a result of the arrival of random information, $W_{k+1}$, consisting of arrivals of customer requests and appearances of crowdshippers. We explain these two transitions in the following two sections. 


\subsection{Actions and Deterministic Transition to Post-Decision State} \label{sec:action}

\noindent At a decision epoch $k$, the dispatcher observes the pre-decision state $s_k$ and selects an action $x \in \mathcal{X}(s_k)$, where $\mathcal{X}(s_k)$ is the set of possible actions corresponding to the pre-decision state $s_k$. An action $x \in \mathcal{X}(s_k)$ consists of two components: (i) it updates the routing plan $\theta_m$ for each vehicle $m \in M_k$ through the assignment and routing of each newly-arriving request $r \in \mathcal{U}_k$ and the possible reassignment and resequencing of outstanding requests $r \in \mathcal{R}_k^o$, and (ii) it specifies the departure time of each vehicle currently idle at a location. For an action to be feasible at epoch $k$, each vehicle $m$ must be able to complete its planned route before the end of its availability, $b_m$, and all requests must be assigned to a vehicle. The presence of one or more dedicated vehicles available until all requests are serviced ensures feasibility. 

The execution of an action $x$ in a pre-decision state $s_k$ incurs a cost $c(s_k, x)$ and triggers a deterministic transition to a post-decision state $s_k^x$. Lines~\ref{line:init}--\ref{line:endwait} of Algorithm~\ref{alg:PDS} break down the execution of action $x$ according to its two components: Line~\ref{line:init} addresses the assignment and sequencing of requests and lines~\ref{line:beginwait}--\ref{line:endwait} determine the temporal scheduling of the next segment of each vehicle's route. Specifically, Line~\ref{line:init} applies the CFA-based DRACE heuristic given state $s_k$ to identify an action $x$ that updates the routing plans for each vehicle through the assignment and routing of outstanding requests. In Section~\ref{sec:DRACE}, we describe DRACE in detail and formally outline it in Algorithm~\ref{alg:action}. 

After the execution of {\tt DRACE} in Line~\ref{line:init}, each request is assigned to a vehicle and each vehicle $m \in \mathcal{M}_k^x$ has an updated routing plan $\theta_m^x$ that maintains vehicle $m$'s current location (or destination if still en route), i.e., $\theta_m^x(1) = \theta_m(1)$. Given the updated routing plans, lines~\ref{line:beginwait}--\ref{line:endwait} then determine when each vehicle $m$ that is currently idle at $\theta_m^x(1)$ will depart for $\theta_m^x(2)$. While updating the timing of the next segment of the route execution, we also calculate the cost $c(s_k,x)$ of the action $x$ applied in state $s_k$. Line~\ref{line:initcost} initializes the $c(s_k,x)$ to zero. 

Line~\ref{line:idlevehicles} identifies the three cases which a vehicle requires updating of its temporal scheduling: (i) if vehicle $m$ has just arrived at $\theta_m(1)$ at epoch $k$ ($f_m = t_k^x$), (ii) if vehicle $m$ is idle at $\theta_m(1)$  (having previously arrived, $f_m < t_k^x$) and action $x$ changed its next location $(\theta_m^x(2) \neq \theta_m(2))$, or (iii) if vehicle $m$ is idle at $\theta_m(1)$  ($f_m < t_k^x$) and the current time is its planned departure time $(w_m = t_k^x)$. 

If the arrival time of a vehicle $m$ at its next destination is beyond the ready time of the request corresponding to the next destination, Lines~\ref{line:loadedLate}--\ref{line:deliverycalc} execute the immediate departure of vehicle $m$ from $\theta_m^x(1)$ and perform the necessary updates. We note this implies that a vehicle will always depart immediately for its next destination if it is a delivery location of a request. When vehicle $m$ departs its current location $\theta_m^x(1)$, Lines~\ref{line:updatearrival} and \ref{line:updatedeparture} set the arrival time and planned departure time for the next destination. Line~\ref{line:addnewtravel} adds the travel cost to the next destination to $c(s_k,x)$. 

Lines~\ref{line:begincheckpickup}--\ref{line:deliverycalc} execute request status updating and cost accounting depending on the next destination's type. If the next destination is a pickup location, then Line~\ref{line:remr} removes the corresponding request from the set of outstanding requests and Line~\ref{line:addr} adds it to the set of in-process requests. If the next destination is a pickup location, lines~\ref{line:begincheckcs} and lines~\ref{line:calcfee} add the per-delivery crowdshipper fee $\rho$ to $c(s_k,x)$ if the servicing vehicle is a crowdshipper. If the next destination is a delivery location, then lines~\ref{line:begincheckdelivery} and \ref{line:deliverycalc} adds any late delivery charge (by multiplying $\mu_2$ by the number of minutes beyond the delivery deadline) to $c(s_k,x)$.

If the location from which a vehicle is immediately departing is a delivery location, then lines~\ref{line:leftdelivery} and \ref{line:remover} remove this request from the set of in-process requests to represent its service completion. Then Line~\ref{line:updateroute} updates $\theta_m(x)$ by removing $\theta_m(1)$, thus making the route one stop shorter.    

If the arrival time of a vehicle $m$ at its next destination is not beyond the ready time of the request corresponding to the next destination, vehicle does not depart immediately and lines~\ref{line:beginwaitcalc} and \ref{line:addwait} set the planned departure time. In this case, the vehicle will plan to wait at its current location at least the amount of time to make its planned arrival time at its next destination coincide with that location's ready time. This minimum amount of wait time reflects the operational restrictions of request pickup. However, we also consider the possibility of additional \emph{strategic} wait time. By having a vehicle wait at a location beyond the time dictated by the request ready time, it may be possible to acquire additional information in the form of new request arrivals and new crowdshipper appearances. This new information may prompt a change to the routing plan that would not have been possible if the vehicle had already committed to the next pickup by departing. We must balance the benefit from the information observed by a vehicle waiting beyond the operational minimum with the opportunity cost of consuming a vehicle's temporal delivery capacity being idle when it alternatively could be progressing through its route and therefore possibly serving future requests in a more timely manner. To determine the planned departure time, Line~\ref{line:addwait} calculates the planned wait time as the larger of: (1) a percentage (given by $\eta$) of the vehicle's remaining time availability, and (2) the operational wait time to arrive at the next destination at its ready time. By calibrating the value of $\eta$ in offline simulation experiments, we consider the balance between too much and too little wait time. 

Finally, Line~\ref{line:stillenroute} shows that the arrival time at and planned departure time from $\theta_m^x(1)$ is left unchanged if vehicle $m$ is still en route ($f_m > t_k^x$) or if vehicle $m$ is idle ($f_m < t_k^x$) and its next destination wasn't changed or its planned departure time hasn't been reached ($w_m > t_k^x$).

\begin{algorithm}[]
\SetAlgoLined
\caption{Dynamics of Markov Decision Process} 
\label{alg:PDS}

\KwData{Pre-decision state $s_k = \left( t_k, \mathcal{R}_k^p, \mathcal{R}_k^o, \mathcal{U}_k, \mathcal{M}_k \right)$ and mechanism to identify an action $x$}
\KwResult{Post-decision state $s_k^x = \left( t_k^x, \mathcal{R}_{k}^{p,x}, \mathcal{R}_{k}^{o,x}, \mathcal{U}_k^x, \mathcal{M}_k^x \right)$}

$s_k^x \leftarrow {\tt DRACE}(s_k)$ \label{line:init}

$c(s_k,x) \leftarrow 0$ \label{line:initcost}

\For{{\bf each }$m \in \mathcal{M}_k^x$} { \label{line:beginwait}
    \uIf{$(f_m == t_k^x)$ {\normalfont or} $\left[ (f_m < t_k^x) \right.$ {\normalfont and} $\left[ (\theta_m^x(2) \neq \theta_m(2)) \right.$ {\normalfont or} $\left.\left. (w_m == t_k^x)\right] \right]$}{ \label{line:idlevehicles}
        \uIf{$t_k^x + \tau(\theta_m^x(1),\theta_m^x(2),t_k^x) \geq  e_{r(\theta_m^x(2))}$}{ \label{line:loadedLate}
            $f_m^x \leftarrow t_k^x + \tau(\theta_m^x(1),\theta_m^x(2),t_k^x)$ \label{line:updatearrival} \\
            $w_m^x \leftarrow f_m^x$ \label{line:updatedeparture} \\
            $c(s_k,x) \leftarrow c(s_k,x) + \mu_1 \tau(\theta_m^x(1),\theta_m^x(2),t_k^x)$ \label{line:addnewtravel}

            \uIf{$\theta_m^x(2) == o_r$ {\normalfont for some} $r \in \mathcal{R}_k^{o,x}$}{ \label{line:begincheckpickup}
                $\mathcal{R}_k^{o,x} \leftarrow \mathcal{R}_k^{o,x} \setminus \left\{ r \right\}$ \label{line:remr} \\
                $\mathcal{R}_k^{p,x} \leftarrow \mathcal{R}_k^{p,x} \cup \left\{ r \right\}$ \label{line:addr} \\
                \uIf{$m \in \mathcal{G}_k^x$}{ \label{line:begincheckcs}
                    $c(s_k,x) \leftarrow c(s_k,x) + \rho$ \label{line:calcfee}       
                }
            }
            
            \uElseIf{$\theta_m^x(2) == d_r$ {\normalfont for some} $r \in \mathcal{R}_k^{p,x}$}{ \label{line:begincheckdelivery}
                $c(s_k,x) \leftarrow c(s_k,x) + \mu_2 \max \left\{ 0, f_m^x - l_r  \right\}$ \label{line:deliverycalc} \\ 
                }
            \uIf{$\theta_m^x(1) == d_r$ {\normalfont for some} $r \in \mathcal{R}_k^{p,x}$}{ \label{line:leftdelivery}
                $\mathcal{R}_k^{p,x} \leftarrow \mathcal{R}_k^{p,x} \setminus \left\{ r \right\}$ \label{line:remover} \\
            }
            Update $\theta_m^x$ by removing $\theta_m^x(1)$ \label{line:updateroute}\\
        }
        \uElse{ \label{line:beginwaitcalc}
            $w_m^x \leftarrow t_k^x + \max \left\{ \eta \left( b_m^x - t_k^x \right), e_{r(\theta_m^x(2))} - t_k^x - \tau(\theta_m^x(1),\theta_m^x(2),t_k^x)\right\}$ \label{line:addwait} \\ 
        } 
    }
    \uElse{
        $f_m^x \leftarrow f_m$, $w_m^x \leftarrow w_m$ \label{line:stillenroute}
    }
}\label{line:endwait}

\vspace{0.5cm}

\KwData{Post-decision state $s_k^x$ and random information $W_{k+1} = \left( \omega_{k+1}, \gamma_{k+1} \right)$}
\KwResult{Post-decision state $s_{k+1}$}

$t_{k+1} \leftarrow t_k^x + 1$ \label{line:incrtime}

$\mathcal{R}_{k+1} \leftarrow \mathcal{R}_k^x$ \label{line:updatereq}

$\mathcal{U}_{k+1} \leftarrow \omega_{k+1}$ \label{line:newreq} 

$\mathcal{M}_{k+1} \leftarrow \mathcal{M}_{k}^x \cup \gamma_{k+1}$  \label{line:newveh}

 \For{{\bf each }$m \in \mathcal{M}_{k+1}$}{ \label{line:timeupdatebegin}
        \uIf{$b_{m} \leq t_{k+1}$ \label{line:expirevehicle} }{
            $\mathcal{M}_{k+1} \leftarrow \mathcal{M}_{k+1} \setminus \left\{ m \right\}$ \label{line:removevehicle}
            } 
        
} \label{line:stochend}

\end{algorithm}

\subsection{Stochastic Transition to Pre-Decision State}

\noindent After the execution of an action $x$ triggering a deterministic transition from pre-decision state $s_k$ to post-decision state $s_k^x$, the arrival of random exogenous information, $W_{k+1}$, triggers a stochastic transition from $s_k^x$ to $s_{k+1}$. In our case, $W_{k+1} = ( \omega_{k+1}, \gamma_{k+1})$, where $\omega_{k+1}$ represents the arrival of a (possibly empty) set of new customer requests and $\gamma_{k+1}$ corresponds to the appearance of a (possibly empty) set of new crowdshippers declaring their availability to perform deliveries.  

Line~\ref{line:incrtime} of Algorithm~\ref{alg:PDS} reflects that we treat the base unit of time to be one minute so that decision epochs occur one minute apart. Line~\ref{line:updatereq} initializes the request information for $s_{k+1}$ from $s_k^x$. Line~\ref{line:newreq} updates the set of newly-arriving requests $\mathcal{U}_{k+1}$ with $\omega_{k+1}$. Line~\ref{line:newveh} supplements the fleet of vehicles $\mathcal{M}_{k+1}$ with $\gamma_{k+1}$. Lines~\ref{line:timeupdatebegin}--\ref{line:stochend} update the status of the fleet by removing any vehicle whose availability has expired.

\section{Solution Methodology}
\label{sec:solution}

\noindent Let $\Pi$ be the set of all Markovian deterministic policies, where a policy $\pi \in \Pi$ is a sequence of decision rules: $\pi = \left( X_0^\pi(s_0), X_1^\pi(s_1), \ldots X_K^\pi(s_K) \right)$ where each decision rule $X_k^\pi(s_k): s_k \rightarrow \mathcal{X}(s_k)$ is a function that specifies the action choice when the process occupies state $s_k$ and follows policy $\pi$, and $K$ represents the epoch index at the end of the problem horizon. We seek a policy $\pi \in \Pi$ that minimizes the total expected cost, conditional on the initial state $s_0$: 
$$
\min_{\pi \in \Pi} \mathbb{E} \left[ \sum_{k=0}^K c(S_k, X_k^\pi(S_k)) | s_0  \right]. 
$$
The Bellman equation recursively restates the cost structure as
$$
V(s_k) = \min_{x \in \mathcal{X}(s_k)} \left\{ c(s_k, x) + \mathbb{E} \left[ V(s_{k+1}) | s_k^x  \right]  \right\} 
$$
for $k = 1, \ldots, K-1$. 

Identifying an optimal policy is challenging due the three curses of dimensionality present in our problem: (1) the state space is multi-dimensional and can grow exponentially, (2) the action space consists of a correspondingly large number of options of assigning and sequencing requests on vehicles, and (3) the outcome space is vast due to the uncertain request arrivals and crowdshipper appearances. Four common strategies for creating policies are: lookahead approximation (LA), value function approximation (VFA), policy function approximation (PFA), and cost function approximation (CFA) (Powell et al. \cite{Powell2012}). LAs can require substantial online computation to assess the future impact of an action, which preclude them from being viable options for real-time problems with large action spaces. VFAs approximate the second term of the Bellman equation (cost-to-go), but are limited by dimensionality, which can grow quite large in routing problems with multiple vehicles. PFAs are best suited when it is possible to design a function (i.e., a policy) that captures the structure of the problem in order to provide decisions. There is no obvious function class to assess the choice of vehicle type (dedicated vehicle or crowdshipper) while also considering the impact of the routing decision on travel cost and lateness charges. Using calibrated parameters, CFAs modify the objective and/or constraints of a deterministic optimization formulation so that it better accounts for future uncertainty.  

We apply a CFA which modifies the first term of the Bellman equation (the current epoch cost) to reflect the future impact of the current action. Given a calibrated parameter $\lambda$, we consider
\begin{equation} \label{eqn:generalcfa}
V(s_k) \approx \min_{x \in \mathcal{X}(_k)} \left\{ \bar{c}_k(s_k, x | \lambda)  \right\}
\end{equation}
for $k = 1, \ldots, K-1$. In the remainder of this section, we describe our solution approach that leverages a CFA to identify a policy.  

\subsection{Destroy-and-Repair Accounting for Capacity Expiration (DRACE)}
\label{sec:DRACE}

\noindent To determine the action $x$ for the system in state $s_k$ at an epoch $k$, we apply a destroy-and-repair heuristic that assigns and sequences requests not associated with any vehicle. In the evaluation of different ways to assign and route the requests, the heuristic uses a CFA to account for the temporal availability of the vehicles. 

Algorithm~\ref{alg:action} outlines the procedure which considers the information in the pre-decision state $s_k$ and implicitly identifies an action $x$ by iteratively assigning and sequencing requests on vehicles. 

Line~\ref{line:initstate} initializes the post-decision state by importing the information from the pre-decision state. In particular, we note that the transition from the pre-decision state to the post-decision state is instantaneous, $t_k^x = t_k$. In Lines~\ref{line:beginrouteupdate}--\ref{line:endrouteupdate}, if a currently scheduled but outstanding request has a ready time within $\Gamma$ time units of the current time, $t_k$, we remove it from its current route. The magnitude of $\Gamma$ affects the required computation time of Algorithm~\ref{alg:action}; we calibrate the value of $\Gamma$ so the per-epoch execution of Algorithm~\ref{alg:action} remains within one minute to facilitate real-time decision-making in our application.

To prioritize the most urgent requests, Line~\ref{line:sortreq} creates an ordered set $\mathcal{I}$ by placing the current unscheduled requests, $\mathcal{U}_k \cup \left\{ \mathcal{R}_k^o: e_r - t_k \leq \Gamma \right\}$, in non-decreasing order of $l_r$. For each request $r \in \mathcal{I}$, Lines~\ref{line:requestsConsidered}--\ref{line:endpfor} iteratively determines a vehicle assignment. For the request $r$ currently being considered, Lines~\ref{line:vehiclesConsidered}--\ref{line:endpfor} estimates the impact of assigning request $r$ to vehicle $m$. Line~\ref{line:deltacalc} executes the {\tt DELTA} procedure to estimate: (1) the increase in vehicle $m$'s travel time due to the assignment of request $r$ ($\Delta_{rm}^{travel}$), (2) the increase in vehicle $m$'s delivery lateness due to the assignment of request $r$ ($\Delta_{rm}^{late}$), and (3) the completion time of vehicle $m$'s route ($t_f$). 

If $r \in \left\{ \mathcal{R}_k^{o}: e_r - t_k \leq \Gamma \right\} \subset \mathcal{I}$, the {\tt DELTA} procedure estimates $\Delta_{rm}^{travel}$, $\Delta_{rm}^{late}$, and $t_f$ via a reconstruction routine that removes all requests from vehicle $m$'s route and then iteratively reinserts request $r$ and the other requests assigned to vehicle $m$. If $r \in \mathcal{U}_k \subset \mathcal{I}$, the {\tt DELTA} procedure estimates $\Delta_{rm}^{travel}$, $\Delta_{rm}^{late}$, and $t_f$ via a fast cheapest insertion routine that searches for the positions to place the pickup and delivery locations for request $r$. We invest less time for the insertion of newly-arriving requests because these requests are generally less urgent than the earlier-arriving outstanding requests. 


Based on the {\tt DELTA} procedure, at epoch $k$, if request $r$ can be added to vehicle $m$ without its route completion time extending beyond its availability window, then lines~\ref{line:vehiclesConsidered}--\ref{line:endfeasible} compute $c_{rm}$, the cost of assigning request $r$ to vehicle $m$, according to the cost function approximation: 
\begin{equation} \label{eqn:equationCFA} 
     c_{rm} = \mu_1 \Delta_{rm}^{travel}+ \mu_2 \Delta_{rm}^{late}+\rho I_{\left\{m \in \mathcal{G}\right\}} + \lambda (b_m - t_k).  
\end{equation}
If vehicle $m$ cannot accommodate request $r$ without its route completion time extending beyond the end of its availability, then we set $c_{rm} = \infty$ in lines~\ref{line:infeasiblebegin} and \ref{line:infeasibleend}. Equation~(\ref{eqn:equationCFA}) consists of four summands capturing the effect of assigning request $r$ to vehicle $m$. The first summand corresponds to the product of the travel time multiplier $\mu_1$ and the increase in vehicle $m$'s travel time due to its assignment of request $r$. The second summand corresponds to the product of the lateness multiplier $\mu_2$ and the increase in vehicle $m$'s violation of request deadlines due to its assignment of request $r$. In the third summand, $I_{\left\{m \in \mathcal{G}\right\}} = 1$ if $m \in \mathcal{G}$ and 0 otherwise implies that the per-delivery fee $\rho$ is charged if $m$ is a crowdshipper and no fee otherwise. The fourth summand corresponds to the product of the CFA parameter $\lambda$ and the amount of time vehicle $m$ is available to service requests. 

The first three summands in Equation~(\ref{eqn:equationCFA}) reflect the economic impact of assigning request $r$ to vehicle $m$. The values of the parameters $\mu_1$, $\mu_2$, and $\rho$ directly come from the economics of the problem context. The value of $\mu_1$ represents the cost of each minute of vehicle travel time. The value of $\mu_2$ represents the cost of each minute of violating a customer's deadline. The value of $\rho$ represents the per-delivery fee remunerating crowdshippers. 

The fourth summand in Equation~(\ref{eqn:equationCFA}) modifies the cost function by introducing a parameterized term that captures the opportunity cost of the limited time availability of the vehicles. While there is not a tangible cost associated with the slack time of a vehicle, accounting for it in this fourth summand rewards actions that assign requests to vehicles whose time availability will expire sooner. This allow greater flexibility in later assignment decisions and also gives the possibility of preserving the temporal capacity of other vehicles in the most critical moments. We determine the CFA parameter $\lambda$ by tuning it to optimize performance in offline simulation experiments. 

In Line~\ref{line:bestVehicle}, we identify vehicle $m^{\star}$ that achieves the smallest value of $c_{rm}$ for request $r$. Line~\ref{line:updateRoute} executes the routing of request $r$ on vehicle $m^{\star}$ by inserting $o_r$ and $d_r$ at their minimum cost positions. Then, the procedure continues with the request $r^{\prime}$ with the next smallest $l_{r^{\prime}}$. 

To relate Algorithm~\ref{alg:action} to Equation~(\ref{eqn:generalcfa}), we first express the cost of an action $x$ inserting the set of unscheduled requests, $\mathcal{U}_k \cup \mathcal{R}_k^{o}$, onto the existing routes of a delivery fleet as  
$$
\bar{c}(s_k, x | \lambda) =  \sum_{m \in \mathcal{M}_k} \left[ \mu_1 \Delta_m^{travel,x} + \mu_2 \Delta_m^{late,x} + \rho q_m^x I_{\left\{m \in \mathcal{G}\right\}} + \lambda q_m^x (b_m - t_k)  \right],   
$$
where $\Delta_m^{travel,x}$ is the increase in vehicle $m$'s travel time due to the requests inserted onto its route by action $x$, $\Delta_m^{late,x}$ is the increase in lateness of vehicle $m$'s deliveries due to the requests inserted onto its route by action $x$, and $q_m^x$ is the number of requests inserted onto vehicle $m$ by action $x$. 

Determining an $\arg\!\min_x \bar{c}_k(s_k, x | \lambda)$ requires solving a deterministic pickup-and-delivery problem which is difficult (NP-hard). To obtain a solution within the one-minute inter-epoch time, DRACE considers a subset of unscheduled requests $\mathcal{U}_k \cup \left\{ \mathcal{R}_k^{o}: e_r - t_k \leq \Gamma \right\}$ and inserts them onto the best vehicle according to Equation~(\ref{eqn:equationCFA}) in increasing order of their delivery deadline. This results in an action $\bar{x}$ with cost 
$$
\bar{c}_k(s_k, \bar{x} | \lambda) = \sum_{r \in \mathcal{I}} \left[ \mu_1 \Delta_{r,m^\star(r)}^{travel} +  \mu_2 \Delta_{r,m^\star(r)}^{late} + \rho I_{\left\{m^{\star}(r) \in \mathcal{G}\right\}} + \lambda (b_{m^\star(r)} - t_k) \right],   
$$
where $m^\star(r)$ is the vehicle on which DRACE inserts request $r$. 
 
\begin{algorithm}[]
\SetAlgoLined
\caption{Destroy-and-Repair Accounting for Capacity Expiration: {\tt DRACE}} 
\label{alg:action}
\KwData{Pre-decision state $s_k$}
\KwResult{Post-decision state $s_k^x$ resulting from an action $x$}

$t_{k}^x \leftarrow t_k$, $\mathcal{R}_{k}^x \leftarrow \mathcal{R}_k$, $\mathcal{U}_{k}^x \leftarrow \mathcal{U}_{k}$, $\mathcal{M}_{k}^x \leftarrow \mathcal{M}_{k} $ \label{line:initstate} 

\For{{\bf each }$m \in \mathcal{M}_k^x$} { \label{line:beginrouteupdate}
    \For{$i = 1$ {\normalfont to} $|\theta_m^x|$}{
        \uIf{$(\theta_m^x(i) == o_r)$ {\normalfont or} $(\theta_m^x(i) == d_r)$ {\normalfont for some} $r \in \left\{ \mathcal{R}_k^{o,x}: e_r - t_k \leq \Gamma \right\}$}{\label{line:beginremovecheck}
            Remove $\theta_m^x(i)$ from $\theta_m^x$ \\ \label{line:removevisit}
        } \label{line:endremovecheck}
    }
} \label{line:endrouteupdate}

$\mathcal{I} \leftarrow \mathcal{U}_k \cup \left\{ \mathcal{R}_k^{o}: e_r - t_k \leq \Gamma \right\}$ ordered in increasing $l_r$ \label{line:sortreq}

\For{{\bf each }$r \in \mathcal{I} $} { \label{line:requestsConsidered}

    \For{{\bf each }$m \in \mathcal{M}_k^x$} { \label{line:vehiclesConsidered}
        $(\Delta_{rm}^{travel}, \Delta_{rm}^{late}, t_{f}) \leftarrow {\tt DELTA}(r, \theta_m^x)$ \\ \label{line:deltacalc}
        
        \uIf{$t_f \leq b_m^x$}{\label{line:feasibilityCheck}
            $c_{rm} \leftarrow \mu_1\Delta_{rm}^{travel}+\mu_2\Delta_{rm}^{late}+\rho I_{\left\{m \in \mathcal{G}\right\}} +\lambda (b_m^x - t_k^x)$ \\ \label{line:cfacalc}
        } \label{line:endfeasible}
        \uElse{ \label{line:infeasiblebegin}
            $c_{rm} \leftarrow \infty$ \\ \label{line:infeasiblecost}
        } \label{line:infeasibleend}
    } \label{line:endmfor}
    $m^\star \leftarrow \arg\!\min\limits_{m \in \mathcal{M}_k^x} \left\{ c_{rm} \right\}$ \\ \label{line:bestVehicle}
    $\theta_{m^\star}^x \leftarrow {\tt INSERT}(r, \theta_{m^\star}^x$)   \label{line:updateRoute} \\
    
}  \label{line:endpfor}
\end{algorithm}

\subsection{Myopic Policy} 

\noindent To benchmark our DRACE approach, we consider a myopic policy based on the pre-decision state application of an adaptive large neighborhood search (ALNS) heuristic. There is a large body of literature demonstrating the effectiveness of ALNS for a wide range of vehicle routing problems. We refer the readers to the recent survey of Mara et al. \cite{Windrasmara2022} for an in--depth overview of the ALNS framework and its applications. In our work, we consider the ALNS of Van der Hagen et al. \cite{van2017pickup} originally designed for a deterministic static pickup and delivery problem and modify it for real-time execution in our stochastic dynamic setting. 


At a decision epoch $k$, the ALNS initializes its solution by considering the set of vehicle routes corresponding to the pre-decision state. 
For each newly-arriving request $r \in \mathcal{U}_{k}$, the ALNS randomly selects a vehicle $m$ for which insertion of $r$ is feasible with respect to the vehicle's remaining time availability, $b_m - t_k$, and inserts $o_r$ and $d_r$ at the end of vehicle $m$'s route, if it is feasible. In its random vehicle selection, the ALNS first considers all crowdshippers and if there is no feasible insertion on a crowdshipper route, then it considers the set of dedicated vehicles. 


To facilitate the real time execution of ALNS, we restrict the set of requests that the ALNS considers in a manner similar to Algorithm~\ref{alg:action}. Specifically, at decision epoch $k$, the ALNS considers the removal and reinsertion of the set of requests that have a ready time with $\Gamma$ time units of the current time, $\left\{ \mathcal{R}_k^o: e_r - t_k \leq \Gamma \right\}$. We calibrate the value of $\Gamma$ and the iteration limit so that the ALNS completes execution within one minute. 

During each iteration, the ALNS removes $n$ requests from the set $\left\{ \mathcal{R}_k^o: \right.$ $\left. e_r - t_k \leq \Gamma \right\}$ and then reinserts these requests. The ALNS considers removing requests through \textit{random removal} and the \textit{Shaw removal}. In random removal, each request has an equal probability to be selected. Shaw removal determines the subset of requests for removal using the concept of \textit{relatedness}. For our problem, the relatedness measure between two customer requests $r_i$ and $r_j$ is: 
\begin{equation} \label{eqn:relatdnessFunc} 
q(r_i, r_j)= \phi (t_{o_io_j}+t_{d_id_j})+\chi (|e_i-e_j|+|l_i-l_j|)
\end{equation}
In (\ref{eqn:relatdnessFunc}), $\phi$ and $\chi$ are constant parameters providing relative weighting of two terms measuring different dimensions of the compatibility of two requests. For a pair of requests $r_i$ and $r_j$, $t_{o_io_j}+t_{d_id_j}$ corresponds to the sum of the average travel time between the pickup locations and the average travel time between the two delivery locations. While in general we consider time-dependent travel times in our problem setting, considering average travel times is sufficient for relatedness because it measures the relative compatibility of two requests' pickup and delivery locations rather than their deadlines. The second term in (\ref{eqn:relatdnessFunc}), $|e_i-e_j|+|l_i-l_j|$, measures the temporal compatibility of the requests' ready times and deadlines. At an ALNS iteration,
Shaw removal selects randomly one of the outstanding requests and then removes $\left\lceil n/2 \right\rceil$ others by identifying them according to the relatedness metric.


After requests have been removed from a solution (either via random removal or Shaw removal), we rebuild the solution using a minimum cost insertion heuristic. The cost calculation considers the sum of routing cost, penalty cost, and crowdshipper fees.
More specifically, considering the removed requests ordered on the basis of the urgency of their delivery deadline, we assess the insertion cost of request on each vehicle. For each vehicle, the insertion heuristic finds the feasible minimum cost insertion position for a request pickup location and then for a request delivery location. To speed-up the computation, we evaluate only the first three positions in each route to re-insert a request pickup location, and only the three positions after the pickup insertion position in which the delivery location can be inserted. We then assign the request to the vehicle with the minimum insertion cost. 

\section{Description of Data Sets} 
\label{sec:instances}
\noindent In this section, we outline the construction of the data sets which we apply our solution methodology. In Section~\ref{sec:requestdata}, we describe the structure of the customer request data in our various instance types. We outline the delivery capacity for our instances, particularly describing crowdshippers in Section~\ref{sec:crowdshippers}. In Section~\ref{sec:traveltimes}, we explain how we model time- and location-dependent travel times in our instances.  

\subsection{Customer Request Data}
\label{sec:requestdata}

\noindent As described by Arslan et al. \cite{arslan2019crowdsourced}, our problem setting corresponds to a \textit{few-to-many} service request pattern in that there are relatively few pickup locations (representing the places of business where pickups are originating) and many delivery locations (representing individual customers). To evaluate our solution approach and the effect of the various characteristics of our problem, we consider four collections of benchmark instances. The first collection of instances, which we refer to as the \emph{Ulmer original}  (UO) instances, comes directly from Ulmer et al. \cite{ulmer2021restaurant} and considers 110 potential pickup (restaurant) locations and 32,000 potential delivery (customer) locations in Iowa City, Iowa. There are three classes of instances in which requests arrive according to homogeneous Poisson process over a seven-hour day at low (25.71 requests/hour), medium (34.29 requests/hour), and high  (42.86 requests/hour) rates, respectively. Requests have a deterministic ready time of 10 minutes after request arrival and a soft deadline of 40 minutes after request arrival. 

Based on the UO instances, we create a second collection of instances, which we refer to as the \emph{nonstationary mixed} (NM) instances, that considers nonstationary customer demand and two types of requests: 50\% short-deadline and 50\% long-deadline. The short-deadline requests are ready 20 minutes after request arrival and have a soft deadline 60 minutes after request arrival. The long-deadline requests are ready 40 minutes after request arrival and have a soft deadline 120 minutes after request arrival. These instances consider same set of 32,000 potential delivery locations as the UO instances, but add 138 new pickup locations to represent long-deadline request origins so that total number of potential pickup locations is 248. There are three classes of instances in which requests arrive according to a nonhomogeneous Poisson process over a ten-hour day at low (22.5 requests/hour), medium (30 requests/hour), and high (37.5 requests/hour) rates, respectively. We provide the details on the hourly arrival rates of the short- and long-deadline requests for these three instance classes in \ref{app:datasets}. 

By modifying the NM instances, we create a third collection of instances, which we refer to as the \emph{$n$-to-1} or \emph{many-to-one} (MTO) instances, that consider the possibility of a single customer placing multiple requests, each with a different pickup location. Multiple requests from a single customer arrive at the same time and all requests in a $n$-to-1 bundle have long-deadlines. We generate the $n$-to-1 Bundle instances by letting each long-deadline customer have approximately a 90\% probability of placing a single request, 7.5\% probability of placing two long-deadline requests, and 2.5\% probability of placing three long-deadline requests. For comparison purposes, we maintain the total number of daily request arrivals and the ratio of long-deadline to short-deadline requests to be the same as in the NM instances. The $n$-to-1 Bundle instances represent a situation in which a customer, either of their own volition or in response to an incentive offered on the e-platform, places requests from multiple businesses which then arrive to the dispatcher at the same time. 

We create a fourth collection of instances, which we refer to as the \emph{1-to-$n$} or \emph{one-to-many} (OTM) instances, we modify the NM data instances to consider the possibility of multiple customers (with delivery locations geographically close to each other) place requests with a single business within a short amount of time. All requests in a 1-to-$n$ arrive within 30 minutes of each other and have long deadlines to delivery locations within 2.414 kilometers of each other. We list the probability of a long-deadline request occurring in a 1-to-$n$ bundle in \ref{app:datasets}. The 1-to-$n$ Bundle instances represent a situation in which a long-deadline request from a customer triggers the circulation of promotions to other customers near the triggering customer's delivery location.

\subsection{Delivery Capacity}
\label{sec:crowdshippers}

\noindent For the NM, OTM, and MTO instances, we consider a delivery workforce of five dedicated vehicles and twenty-eight crowdshippers available during randomly-generated one- to four-hour time windows. We generate the appearance times of the twenty-eight crowdshippers according to a homogeneous Poisson process with a rate of three crowdshipper appearances per hour. Each crowdshipper appears to the system at a randomly-generated origin and destination (specifying where the crowdshipper needs to be located at the end of their availability time window). This destination prevents a request assignment that sends a crowdshipper to a distant delivery location near the end of their time window. We provide the availability information for the set of twenty-eight crowdshippers we consider in \ref{app:datasets}. 

For the UO instances, we consider a delivery workforce consisting of three dedicated vehicles available over the entire time horizon and twenty-two crowdshippers available in randomly-generated one- to four-hour time windows. We specify this delivery capacity to be similar to the number of hours of delivery capacity considered by Ulmer et al. \cite{ulmer2021restaurant} which considers a fleet of 15 dedicated vehicles. The availability information for these crowdshippers is the first 22 entries listed in \ref{app:datasets}.

\subsection{Time-Dependent Travel Time}
\label{sec:traveltimes}

\begin{figure}[htb]
\caption{Eight speed profile regions of Iowa City.}
\label{figura:IowaAreas}
\begin{center}
\includegraphics[height=2.5in,width=4in,angle=0]{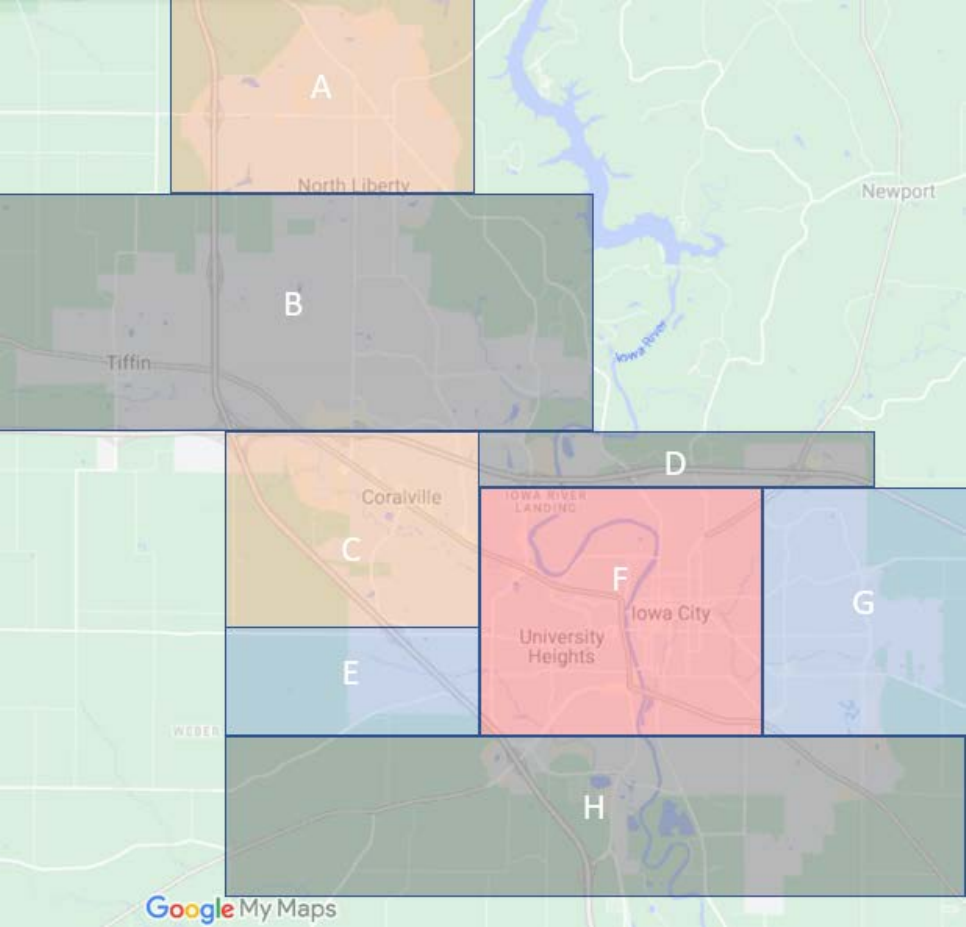}
\end{center}
\end{figure}

\noindent Traffic congestion depends on location and time of day. When possible, the effects of traffic congestion may be mitigated by avoiding being at the wrong place at the wrong time. One strategy to achieve this is to incorporate time--dependent travel in the routing model. We implement the time-dependent nature of travel using a speed profile as in Sun et al.\cite{sun2018time}, and extend it to also vary by geographic location so that $v_{rt}$ corresponds to the travel speed of a vehicle within region $r$ at time $t$. For our instances, Table~\ref{table:speedprofile} lists the 32 travel speed values for $v_{rt}$ corresponding to the four different time periods and eight different geographic regions of the Iowa City area. Figure~\ref{figura:IowaAreas} provides a map displaying the eight speed profile regions of the Iowa City area. We provide further details of our travel-time implementation in \ref{app:datasets}. 


\begin{table}[tbh]
\caption{Speed profile of $v_{rt}$ values (km/min).}
\label{table:speedprofile}
\centering
\begin{tabular}{ccccc } 
& \multicolumn{4}{c}{Time Period ($w$)} \\
\cline{2-5}
Region ($r$) & 8:00 - 10:00 & 10:00 - 18:00 & 18:00 - 20:00 &  $>$ 20:00\\
\hline
 A &  0.25 & 0.40 & 0.25 & 0.40\\ 
 B &  0.50 & 0.67& 0.50& 0.67 \\ 
 C &  0.25 & 0.40 & 0.25 & 0.40\\
 D &  0.50 & 0.67& 0.50& 0.67\\ 
 E &  0.33 & 0.53 & 0.33 & 0.53\\ 
 F &  0.16 & 0.26 & 0.16 & 0.26 \\ 
 G &  0.33 & 0.53 & 0.33 & 0.53\\ 
 H &  0.50 & 0.67& 0.50& 0.67\\ 
 \hline
\end{tabular}
\end{table}

\section{Computational Results}
\label{sec:computation}

\noindent There are four objectives of our computational study. First, we demonstrate the effectiveness of DRACE relative to a myopic ALNS approach. Second, we illustrate the impact of allowing vehicles to strategically wait at locations as executed in lines~\ref{line:beginwait}--\ref{line:endwait} of Algorithm~\ref{alg:PDS} by comparing to vehicle scheduling with only the minimum operational waiting time. Third, we examine the effect of different demand patterns by comparing solution quality for the nonstationary mixed (NM), the many-to-one (MTO), and one-to-many (OTM) instances. Finally, we investigate the impact on solution quality of modeling time-dependent travel-time versus using average travel time. 

In these computational results, we tune parameter values $\lambda = 0.05$ in Line~\ref{line:cfacalc} of Algorithm~\ref{alg:action} and $\eta = 0.20$ in Line~\ref{line:addwait} of Algorithm~\ref{alg:PDS} via offline simulations. We set the cost parameters as $\mu_1 = 1$, $\mu_2 = 5$, and $\rho = 2$. For the myopic ALNS, we calibrate $\phi = 9$ and $\chi = 3$ in Equation~(\ref{eqn:relatdnessFunc}) via offline simulations. To complement the summaries we describe in this section, we provide extended computational results in \ref{app:moreresults}. 


\subsection{DRACE versus Myopic ALNS}
\label{secsec:intraInstancesComp}

\noindent For the four different instance types (NM, OTM, MTO, UO) and three different demand levels (low, medium, high), we conduct computational tests on 100 different daily instances. We compare DRACE to the myopic ALNS with respect to the average total cost (the sum of travel cost, late penalty, and crowdshipper fees). In addition, we compare the two approaches with respect to average minutes late per request to observe how they perform with respect to timely service. As a third measure of comparison, we also compute the percentage of requests served by crowdshippers for both approaches.

As Table~\ref{table:avgcostdelay} displays, DRACE outperforms the myopic ALNS with respect to average total cost and average lateness per request for all type/demand settings. While the average total cost and average lateness per request causes an increase for both approaches, the increase in cost for DRACE is much smaller than for the myopic ALNS. 

A primary difference between DRACE and the myopic ALNS is DRACE's differentiation of delivery vehicles by their remaining temporal delivery capacity via the last term in Equation~(\ref{eqn:equationCFA}) when assigning requests to vehicles. Comparing the percentage of requests served by crowdshippers in the last two columns of Table~\ref{table:avgcostdelay}, we observe the CFA mechanism of Equation~(\ref{eqn:equationCFA}) guides DRACE to utilize the crowdshippers more than the myopic approach for all type/demand settings. For the NM instances, DRACE assigns 94\% of requests to crowdshippers in the low demand case, 80\% of requests to crowdshippers in the medium demand case, and 77\% of requests to crowdshippers in the high demand case, with a similar pattern holding for the OTM and MTO instances. As the number of total requests increases with an increase in demand rate, DRACE inevitably assigns a smaller percentage of requests to crowdshippers because dedicated vehicles must serve more requests as the crowdshippers' capacity approaches its maximum utilization. However, the myopic ALNS assigns roughly the same percentage of requests to crowdshippers as demand rate increases because the crowdshippers still have ample capacity under the myopic approach's request assignment. For the UO instances, DRACE's percentage of requests served by crowdshippers also decreases as the demand rate increases although more gradually due to the shorter deadlines and slightly higher average demand rate of the low demand case (homogeneous rate of 25.71 requests/hour versus heterogeneous rate with an average of 22.5 requests/hr). For the UO instances, the myopic ALNS assigns a larger percentage of requests to crowdshippers as the demand rate increases to stem the lateness fees resulting in this capacity-stressed system. 


\begin{table}
\caption{Comparison of DRACE and myopic ALNS (My.) for four instance types and three demand levels.}
\label{table:avgcostdelay}
\centering
\begin{tabular}{l l r r r r r r } 
 \hline
 & & \multicolumn{2}{c}{Avg. Total} & \multicolumn{2}{c}{Avg. Lateness} & \multicolumn{2}{c}{\% Requests}\\
 & & \multicolumn{2}{c}{Cost/Request} & \multicolumn{2}{c}{Per Request} & \multicolumn{2}{c}{Crowdserved}\\
\hline
Instance & Demand & My. & DRACE & My. & DRACE & My. & DRACE \\
\hline
      &  low & 19 & 13 & 8 & 4 & 40 & 94 \\ 
 NM &  medium & 40 & 15 & 30 & 7& 38 & 80 \\ 
      &  high & 74 & 16 & 51 & 8 & 38 & 77 \\
  \hline
     &  low & 21 & 13 & 9 & 4 & 40 & 93 \\ 
 OTM &  medium & 36 & 18 & 27 & 8 & 39 & 80 \\ 
     &  high & 69 & 20 & 49 & 8 & 38 & 77 \\ 
   \hline
     &  low & 18 & 12 & 8 & 4 & 39 & 93 \\ 
 MTO &  medium & 37 & 15 & 27 & 7 & 38 & 81 \\ 
     &  high & 71 & 16 & 49 & 8 & 38 & 78 \\ 
  \hline
     &  low & 81 & 16 & 41 & 11 & 39 & 90 \\ 
 UO & medium & 163 & 26 & 59 & 19 & 47 & 87 \\ 
     &  high & 230 & 33 & 69 & 22 & 53 & 85 \\  
 \hline
\end{tabular}
\end{table}

\begin{figure}[]
\caption{Percent reduction in average cost per request by DRACE versus myopic ALNS.}
\label{figura:diffavgtotalcost}
\begin{center}
\includegraphics[height=6.25cm,width=6.25cm,angle=0]{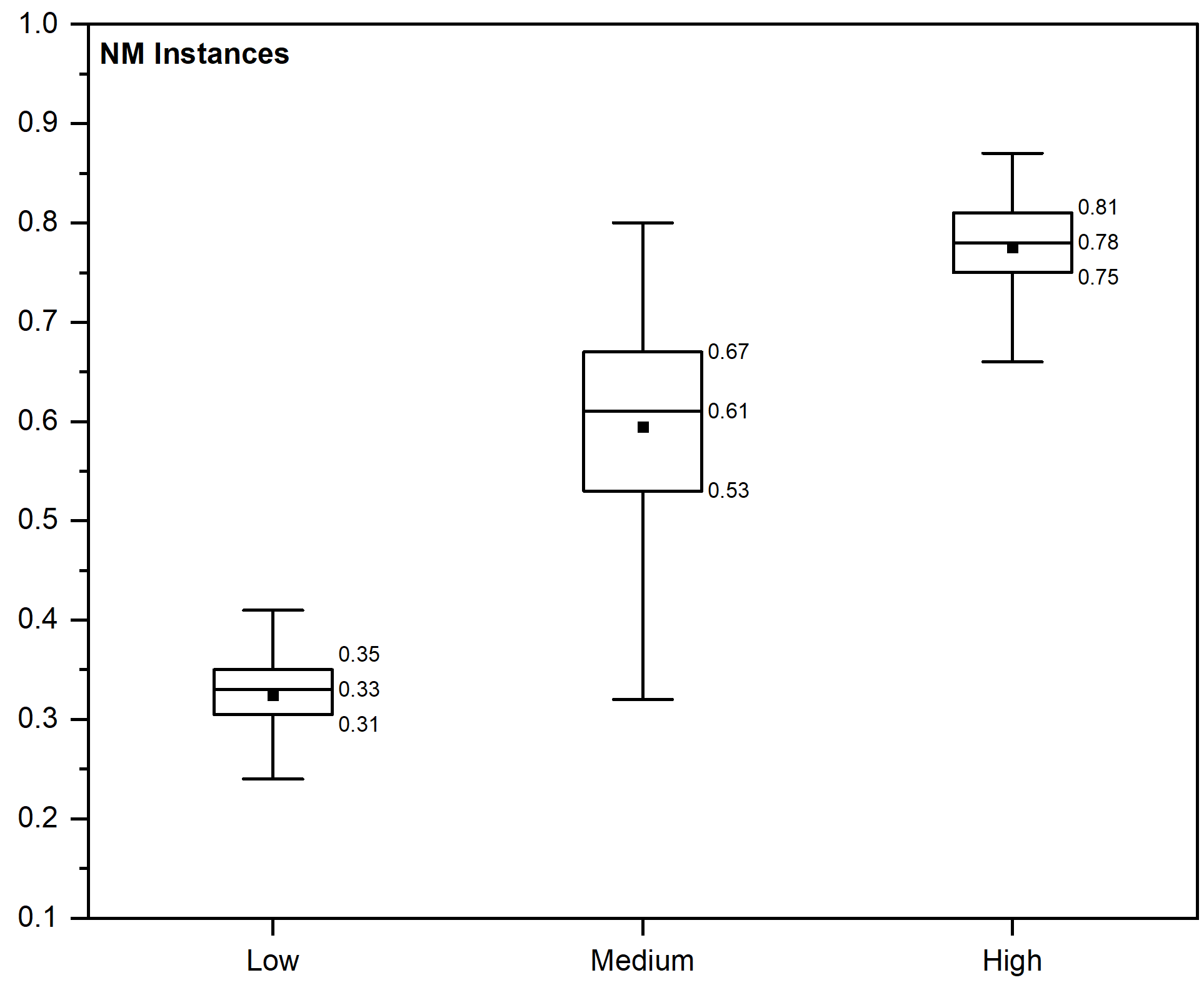}
\includegraphics[height=6.25cm,width=6.25cm,angle=0]{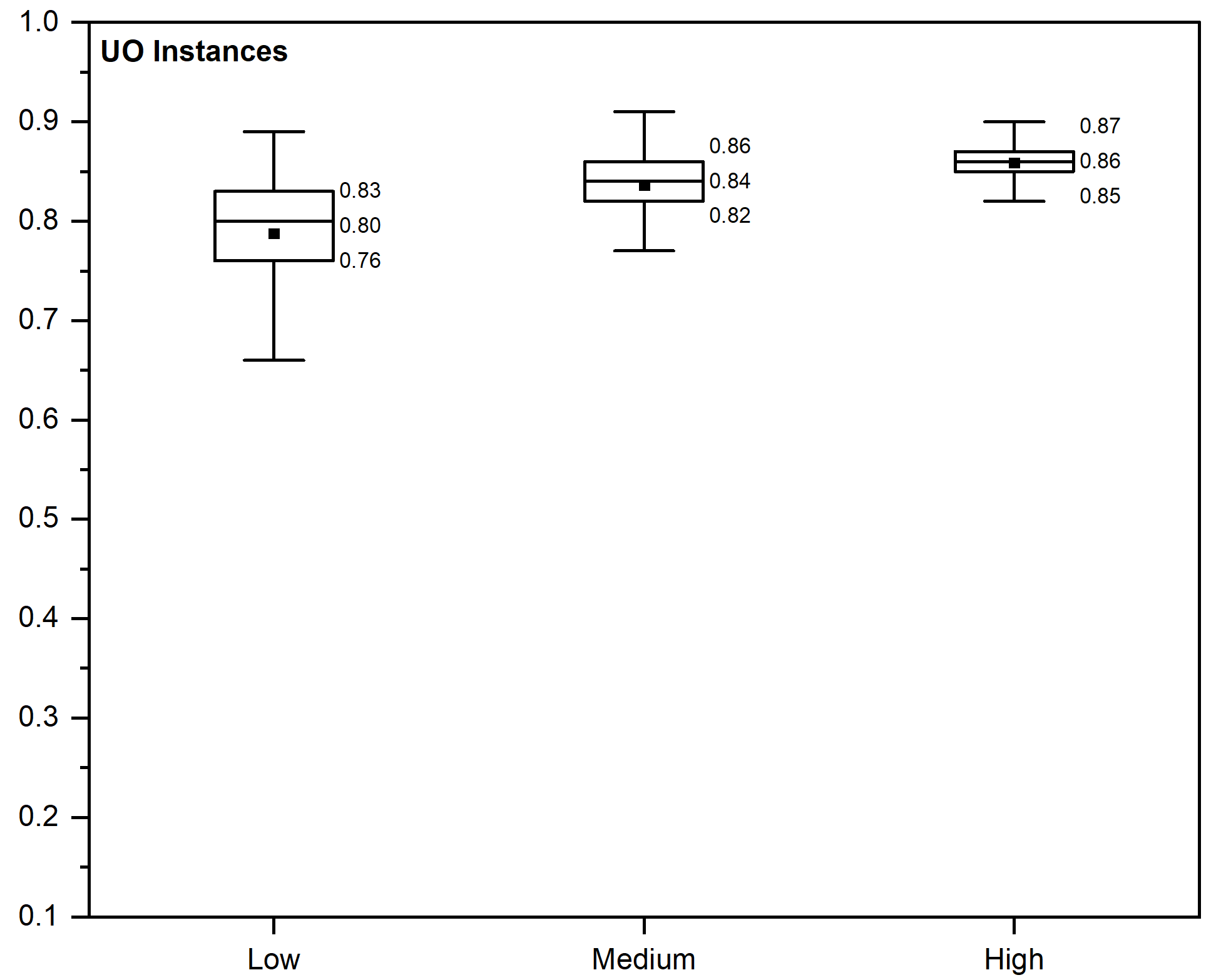}
\end{center}
\end{figure}

Figure~\ref{figura:diffavgtotalcost} displays box plots of the pairwise percent reduction in average cost per request achieved by DRACE versus the myopic approach as computed by \\
$\left( cost_{myopic} - cost_{DRACE} \right)/ cost_{myopic}$. For the NM instances, DRACE achieves a median percent reduction of 33\% in the low demand case, 61\% in the medium demand case, and 78\% in the high demand case, with a similar pattern holding for the OTM and MTO instances (see \ref{app:morefigures}). The minimum percent reduction values are positive for all type/demand settings, demonstrating that DRACE achieves a lower total cost in all 400 instances. The shorter delivery deadlines (30 minutes after request ready time and 40 minutes after request arrival) of the UO instances amplify DRACE's cost advantage over the myopic ALNS. For the UO instances, DRACE achieves a median percent reduction of 80\% in the low demand case, 84\% in the medium demand case, and 86\% in the high demand case. Examining the variation in the distribution of the percent reduction values for the UO instances, we observe that as demand rate increases, DRACE's cost advantage becomes more consistent.

\begin{figure}[]
\caption{Percent reduction in average lateness per request by DRACE versus myopic ALNS. }   
\label{figura:diffAvgDelay1}
\begin{center}
\includegraphics[height=6.25cm,width=6.25cm,angle=0]{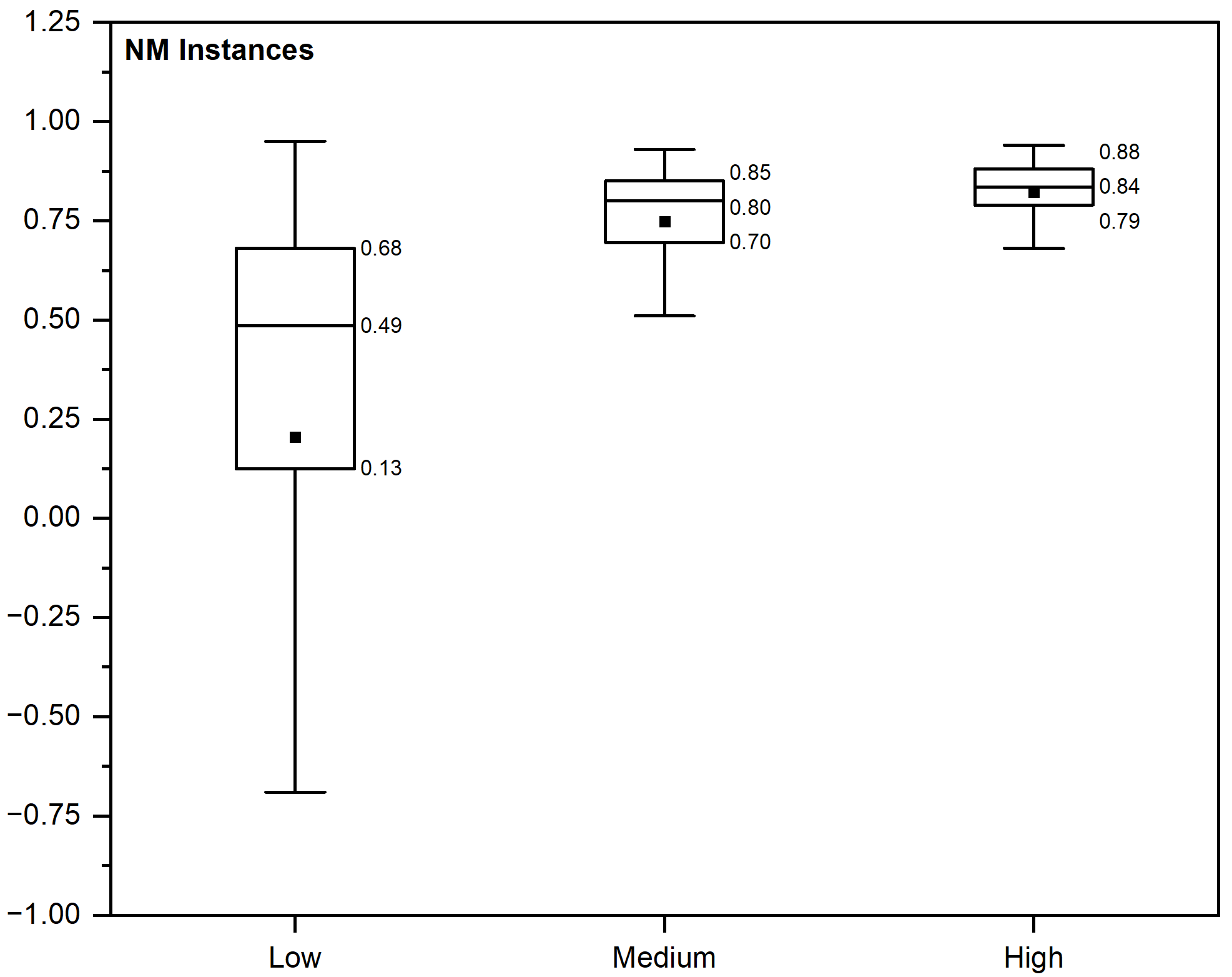}
\includegraphics[height=6.25cm,width=6.25cm,angle=0]{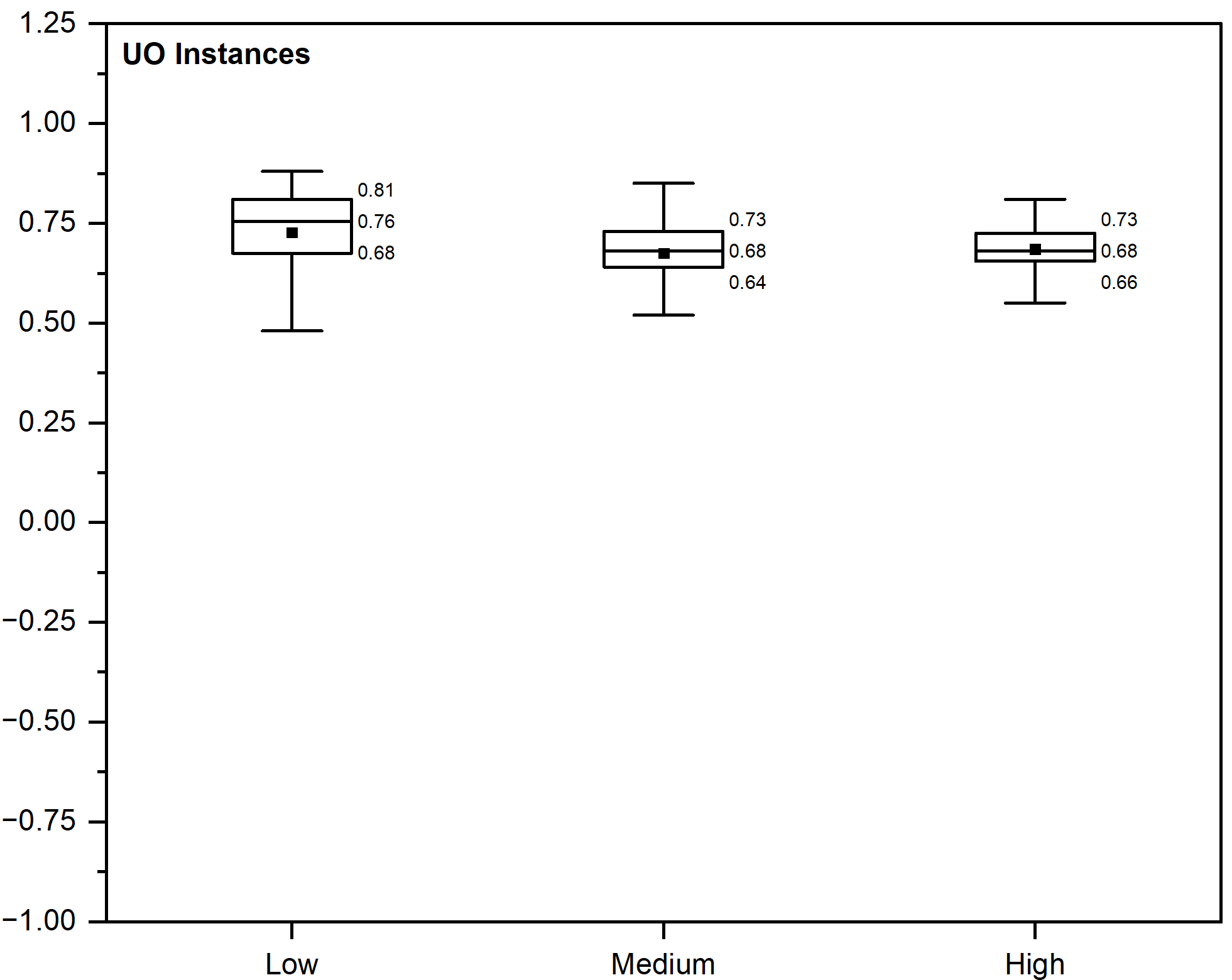}
\end{center}
\end{figure}

Figure~\ref{figura:diffAvgDelay1} displays box plots of the pairwise percent reduction in average lateness per request achieved by DRACE versus the myopic approach as computed by \\$\left( lateness_{myopic} - lateness_{DRACE} \right)/ lateness_{myopic}$. For the NM instances, DRACE achieves a median percent reduction of 49\% in the low demand case, 80\% in the medium demand case, and 84\% in the high demand case, with a similar pattern holding for the OTM and MTO instances (see \ref{app:morefigures}). For the low demand case of the NM, OTM, and MTO instances, there are cases in which the myopic approach results in a smaller average lateness per request than DRACE. In these rare occasions, the CFA mechanism in Equation~(\ref{eqn:equationCFA}) guides DRACE to utilize crowdshipper capacity even though it results in more lateness and the conserved dedicated vehicle capacity is not worth this tradeoff due to the low demand rate. However, as the demand rate increases in the NM, OTM, and MTO instances, there are no cases of DRACE resulting in more lateness than the myopic ALNS as the increased utilization of crowdshipper capacity becomes progressively more important. The shorter delivery deadlines (30 minutes after request ready time and 40 minutes after request arrival) of the UO instances amplify DRACE's service advantage over the myopic ALNS. For the UO instances, DRACE achieves a median percent reduction of 76\% in the low demand case, 68\% in the medium demand case, and 68\% in the high demand case. For the UO instances, as demand rate increases DRACE's service advantage over the myopic ALNS decreases slightly and levels off. This occurs because the shorter deadlines already result in a capacity-stressed system in the low demand case and further increasing demand does not result in additional opportunities for DRACE's higher utilization of crowdshippers and conserved dedicated vehicle capacity to further increase the service gap between DRACE and the myopic ALNS.

\subsection{Impact of Strategic Waiting}

\noindent As we describe in Section~\ref{sec:action}, there may be a benefit for a vehicle to wait at a location and gain new information in the form of new request arrivals and new crowdshipper appearances. By waiting rather than departing for a request pickup, it may be be possible to modify a vehicle's routing plan reduce the total cost. Alternatively, waiting too long can waste a vehicle's temporal delivery capacity and impair its ability to serve requests in a timely manner.

\begin{figure}[]
\caption{Percent reduction in average cost per request with strategic waiting versus no strategic waiting. }
\label{figura:diffwt_cfa_mto}
\begin{center}
\includegraphics[height=6.25cm,width=6.25cm,angle=0]{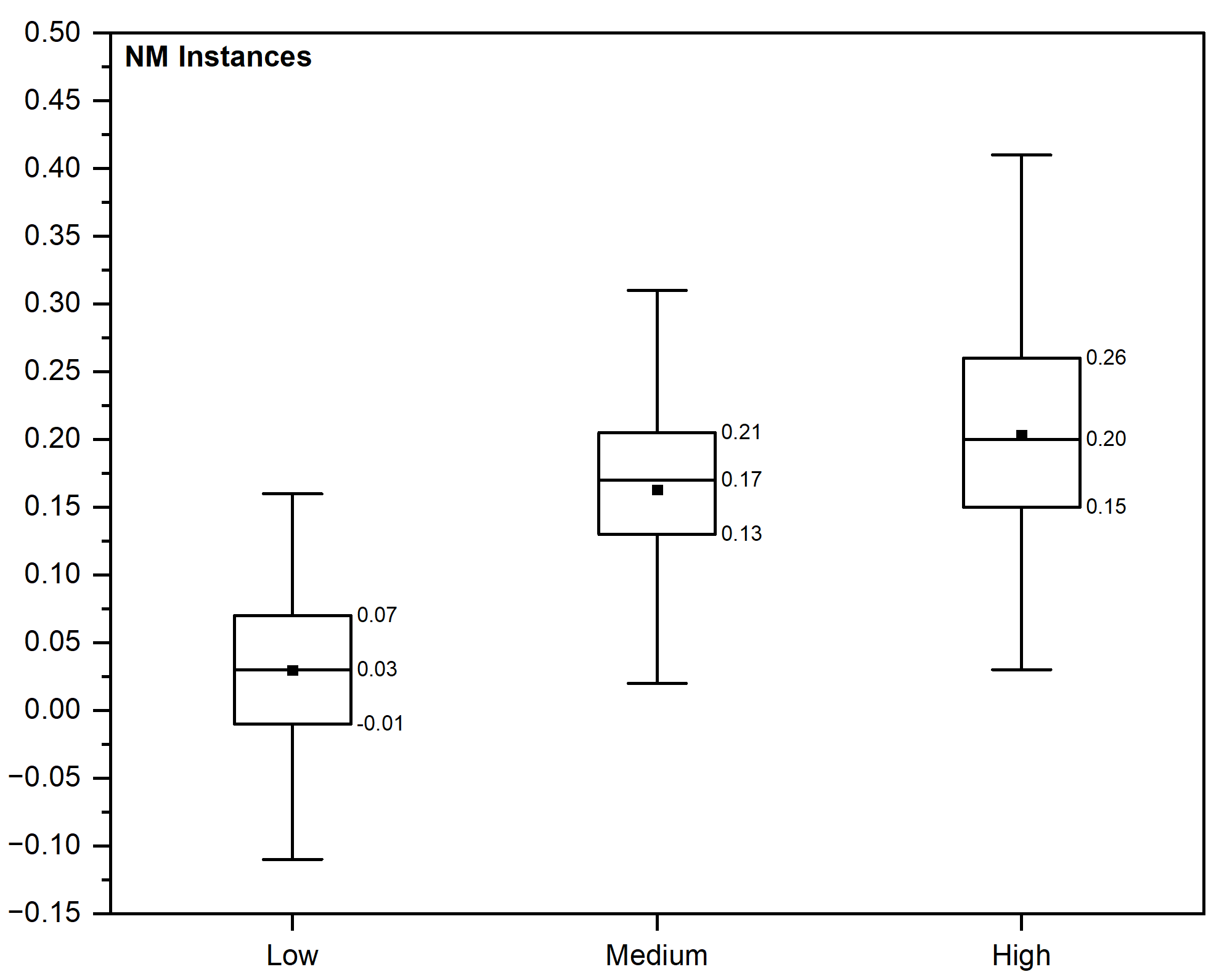}
\end{center}



\end{figure}

Governed by the value of the $\eta$ parameter, lines~\ref{line:beginwaitcalc} and \ref{line:addwait} of Algorithm~\ref{alg:PDS} consider the possibility of allowing a vehicle to strategically wait at a location beyond the operational minimum. In this section, we examine the total cost of solution with strategic waiting corresponding to the tuned value of $\eta = 0.20$ versus the total cost of a solution with only operational waiting and no strategic waiting ($\eta = 0.00$). For the NM instances, Figure~\ref{figura:diffwt_cfa_mto} displays the percent reduction in average cost per request when executing DRACE with $\eta = 0.20$ versus $\eta = 0.00$. The median percent reduction in average cost per request resulting from strategic waiting is 3\% for the low demand case, 17\% for the medium demand case, and 20\% for the high demand case. As the demand rate increases, strategic waiting is increasingly beneficial because more information (more new request arrivals) is observed. For the low demand cases, there are several instances in which the routing with strategic waiting is more costly than routing without strategic waiting as the new information observed during the additional wait time does not merit the time spent waiting to acquire it.   


\subsection{Impact of Demand Management Mechanisms}

\noindent In this section, we analyze the difference in the performance of DRACE across the NM, MTO, and OTM instance types. Recall that that MTO and OTM instances modify the NM instances via two different demand management mechanisms. For the MTO instances, ``many-to-one" requests consisting of a single customer placing requests with multiple businesses are more prevalent than in the base MN instances; approximately 5\% of total requests arrive with one or more other requests from the same customer. For the OTM instances, ``one-to-many" requests consisting of multiple different customers (co-located within specified radius) placing requests with the same pickup location within a thirty-minute period; approximately 11\% of total requests arrive within one of these 1-to-$n$ bundles. 

  

Table~\ref{table:avgcostdelayhetreq} contains 95\% confidence intervals on the pairwise differences in average total cost between the NM instances and MTO instances, and between the NM instances and MTO instances for all three demand levels. For the low demand case, the both the $n$-to-1 and 1-to-$n$ bundling results in a statistically significant reduction in the average total cost. The $n$-to-1 bundling results in reduction between 0.8\% to 5.4\% of the average total cost in the low demand NM instances. The 1-to-$n$ bundling results in reduction between 0.04\% to 5.55\% of the average total cost in the low demand NM instances. 

As the demand rate increases, Table~\ref{table:avgcostdelayhetreq} shows that the increasingly stressed delivery capacity is unable to consistently leverage the $n$-to-1 and 1-to-$n$ bundling, thereby eliminating any benefit of these demand management mechanisms. Therefore, we conclude that these demand management mechanisms are only effective at lowering costs when delivery capacity is sufficiently sized to take advantage of the tested density of bundled requests. We also note that increasing the density of bundled requests should decrease the delivery capacity required to exploit this request arrival pattern.



\begin{table}[]
\centering
\caption{95\% confidence intervals on the effect of $n$-to-$1$ and $1$-to-$n$ request bundling on total daily cost.}
\label{table:avgcostdelayhetreq}
\begin{tabular}{l r r r } 
 \hline
 &  \multicolumn{3}{c}{Demand}\\
\hline
  Instance Comparison     & Low & Medium & High\\
\hline
  NM - MTO  &  [24, 151] & [-120, 124] & [-92, 334]\\ 
  NM - OTM   & [1, 157] & [-152, 84] & [-110, 328]  \\ 
  \hline
  
 \hline
\end{tabular}
\end{table}

%

\subsection{Impact of Modeling Time-Dependent Travel Time}

\begin{figure}[]
\caption{Comparing costs of policies using time-dependent travel time versus cost of policies using average travel time.}
\label{figura:td_check}
\begin{center}
\includegraphics[height=7cm,width=8cm,angle=0]{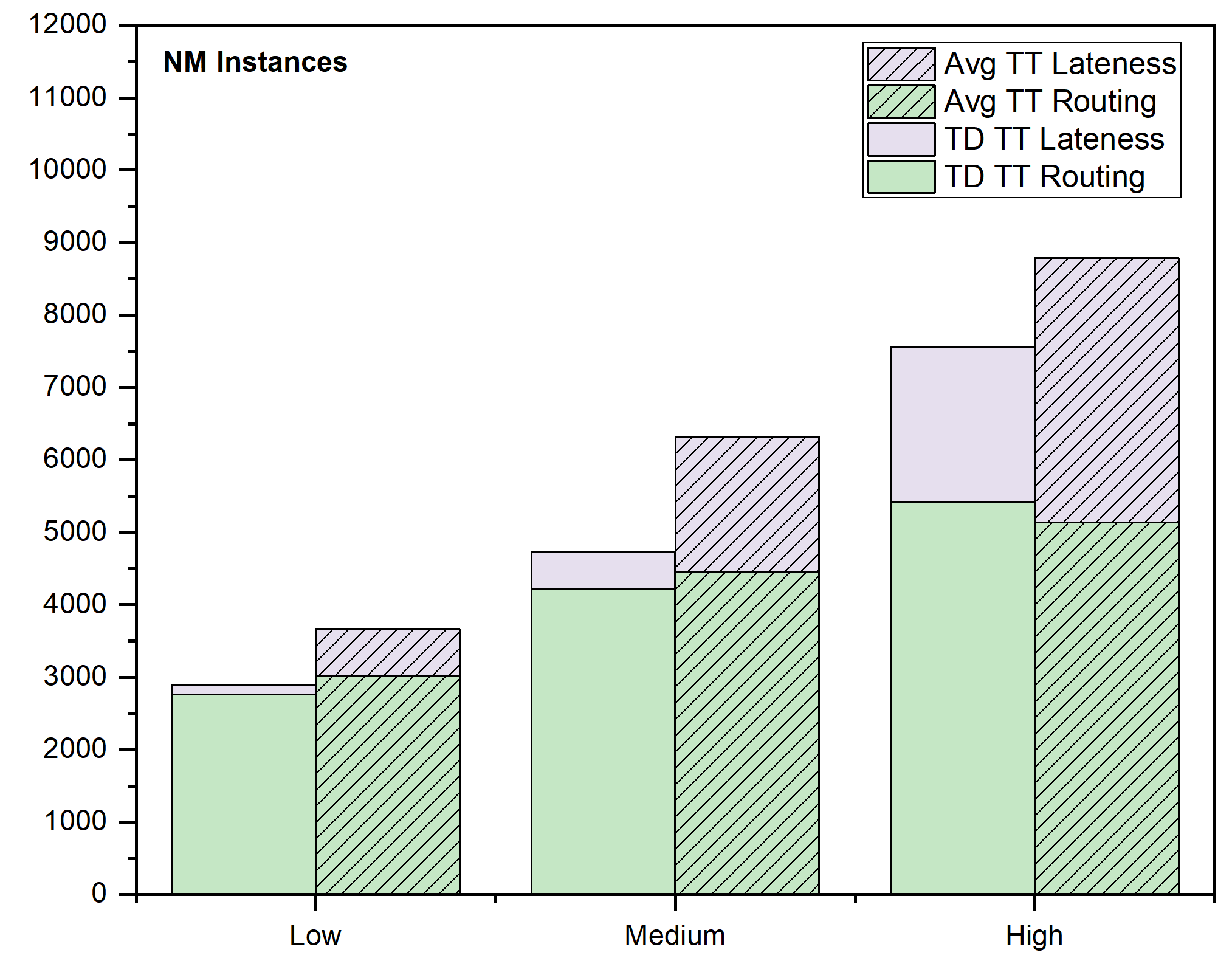}

\end{center}
\end{figure}

\noindent In this section, we examine the benefit of explicitly accounting for the time-dependent travel time in our instances (described in Section~\ref{sec:traveltimes}). We compare the cost of the routing policies that incorporate time-dependent travel time in Table~\ref{table:speedprofile} when determining an action for a current state to the cost of the routing policies that assume an average travel time of $0.4\bar{3}$ km/min when determining an action for a current state. For the NM instances, Figure~\ref{figura:td_check} plots the cost of both the time-dependent travel time (TD TT) policies and the average travel time (Avg TT) policies. Figure~\ref{figura:td_check} breaks down the total cost into two components: the routing cost (which includes fees paid to crowdshippers) and the lateness charge. As Figure~\ref{figura:td_check} shows, policies that consider time-dependent travel times achieve smaller total costs than their average travel time counterparts for all type/demand settings. Policies that consider time-dependent travel times reduce total cost by an average of 21\% over all type/demand settings. The reduction in total cost by considering time-dependent travel times relies primarily on the reducing lateness charge versus reducing routing cost. TD TT policies reduce lateness charges by an average of 63\% but only reduce routing costs by an average of 5\% over all type/demand settings relative to Avg TT policies. 

For the high demand NM instances, considering time-dependent travel time actually increases the average routing cost by 6\% compared to the Avg TT policies. However, the TD TT policies reduce the average lateness charge by 42\% in the high demand NM instances, resulting in a decrease of 14\% in total cost relative to the Avg TT policies. 

Overall, the results of Figure~\ref{figura:td_check} demonstrate that accurately capturing time-dependent travel times in our dynamic pickup-and-delivery problem plays a key role in efficiently assigning requests to vehicles. This finding is important given the computational complications required to account for time-dependent travel times, particularly in a real time setting such in which decisions must be executed on a minute-by-minute basis.


\section{Conclusion}
\label{sec:concl}

\noindent In this study, we consider a dynamic pickup-and-delivery problem in which delivery capacity is supplied through a fleet of dedicated vehicles and crowdshippers. The appearance time of crowdshippers is unknown a priori, but their availability duration is declared upon appearance. To maintain practical relevance, we incorporate time-dependent and location-dependent travel times to execute request assignment and scheduling actions within a one-minute inter-epoch time. 

To form our real time policies to assign and schedule requests, we design a destroy-and-repair heuristic (DRACE) that accounts for the capacity expiration of the delivery vehicles. Our computational results demonstrate that DRACE achieves a substantial percent reduction in total cost relative to a myopic procedure and that this advantage increases as the demand rate increases. Further, we show the impact of inserting strategic wait time into the scheduling of the pickup-and-delivery routes and that the benefit of our wait time mechanism increases as the demand rate increases. We also investigate the impact of two demand management mechanisms that facilitate many-to-one and one-to-many request bundles, respectively. Our computational results demonstrate that both of these demand management mechanisms can improve the efficiency of servicing requests, but only when there exists sufficient delivery capacity relative to the demand rate. Finally, we justify modeling time-dependent travel times by showing that this explicit modeling reduces the total cost by 21\% over all our instances relative to modeling average travel times. 

Our model and solution approach are general enough to accommodate various types of committed and uncommitted crowdshippers Savelsbergh \& Ulmer\cite{Savelsbergh2022} including: gig workers committing to availability time windows known a priori to the dispatcher, (ii) in-store customers who declare themselves as eligible to deliver packages of other customers \cite{dayarian2020crowdshipping}, and (iii) commuters \cite{arslan2019crowdsourced}. Future research remains to compare and contrast these various types of crowdshippers in a workforce scheduling context \cite{ulmer2020workforce}.



\renewcommand{\baselinestretch}{1.0}
\normalsize

\bibliographystyle{plain}

\newpage

\appendix

\section{Data Set Details}
\label{app:datasets}

\noindent In this appendix, we provided details regarding the data sets upon which we apply our solution approach. 

In the NM instances, which consider a request mix of 50\% short deadlines and 50\% long deadlines, the arrival rate of short-deadline requests is non-stationary over the day while the arrival rate of the long-deadline requests is stationary. This captures the difference in demand pattern between requests for perishable items that commonly compose short-deadline requests and requests for non-perishable items that compose long-deadline requests. Table~\ref{table:nonhPP} specifies the hourly arrival rates of the short- and long-deadline requests for these three demand classes.

\begin{table}[htp]
\caption{Hourly rates for nonhomogenous Poisson process in NM instances.}
\label{table:nonhPP}
\centering
\begin{tabular}{r r r r r r r} 
& \multicolumn{6}{c}{Demand Class} \\
\cline{2-7}
     & \multicolumn{2}{c}{Low} & \multicolumn{2}{c}{Medium} & \multicolumn{2}{c}{High} \\
 \hline
Hour & Short & Long & Short & Long & Short & Long \\
 1 & 3.75 & 11.25 & 5 & 15 & 6.25 & 18.75 \\ 
 2 & 11.25 & 11.25  & 15 & 15 & 18.75 & 18.75 \\ 
 3 & 18.75 & 11.25 & 25 & 15 & 31.25 & 18.75 \\ 
 4 & 15.00 & 11.25 & 20 & 15 & 25.00 & 18.75 \\ 
 5 & 11.25 & 11.25 & 15 & 15 & 18.75 & 18.75 \\ 
 6 & 3.75 & 11.25 & 5 & 15 & 6.25 & 18.75 \\ 
 7 & 3.75 & 11.25 & 5 & 15 & 6.25 & 18.75 \\ 
 8 & 11.25 & 11.25 & 15 & 15 & 18.75 & 18.75 \\ 
 9 & 18.75 & 11.25 & 25 & 15 & 31.25 & 18.75 \\ 
 10 & 15.00 & 11.25 & 20 & 15 & 25.00 & 18.75 \\ 
 \hline
\end{tabular}
\end{table}

For the 1-to-$n$ Bundle instances, Table~\ref{table:1ton} lists the probability of long-deadline 1-to-$n$ bundles occurring for values of $n$ from 1 to 6. Both the yield percentage of a promotion and the customer density affect the distribution of 1-to-$n$ bundles. Our instances correspond to neighborhoods of 70 customers with a yield percentage of 0.35\%. 

\begin{table}
\caption{Probabilities of a long-deadline request occurring in a 1-to-$n$ bundles.}
\label{table:1ton}
\centering
\begin{tabular}{r r } 
$n$ & Probability  \\
\hline
 1 & 0.782368  \\ 
 2 & 0.192354 \\ 
 3 & 0.023308  \\ 
 4 & 0.001856 \\ 
 5 & 0.000109 \\ 
 6 & 0.000005\\ 
 \hline
\end{tabular}
\end{table}

\subsection{Delivery Capacity}

\noindent Table~\ref{table:crowddata} contains information on the crowdshipper appearance pattern we use in our computational experiments. 
{\relsize{-2}
\begin{table}
\caption{Crowdshipper $g$ availability information.}
\label{table:crowddata}
\centering
\begin{tabular}{rrrrrrr} 
& &  & \multicolumn{2}{c}{Origin} & \multicolumn{2}{c}{Destination} \\
\cline{4-5} \cline{6-7}
$g$ & $a_g$ & $b_g$ &		Latitude	& Longitude	&	  Latitude	&	Longitude \\	
1	&	1	&	120	&	41.63562541	&	-91.51196350	&	41.65909800	&	-91.55525400	\\
2	&	60	&	180	&	41.63562541	&	-91.51196354	&	41.69834515	&	-91.58892941	\\
3	&	70	&	310	&	41.64128623	&	-91.56508335	&	41.65806850	&	-91.53102480	\\
4	&	90	&	270	&	41.64885944	&	-91.55320966	&	41.64506561	&	-91.52355511	\\
5	&	115	&	295	&	41.71443676	&	-91.58289955	&	41.66152363	&	-91.47081150	\\
6	&	115	&	355	&	41.67312935	&	-91.57551051	&	41.72164463	&	-91.59286545	\\
7	&	130	&	250	&	41.65001141	&	-91.46857959	&	41.72164463	&	-91.59286545	\\
8	&	135	&	315	&	41.63917026	&	-91.51290389	&	41.76164463	&	-91.69286545	\\
9	&	140	&	320	&	41.66699952	&	-91.48110701	&	41.63562541	&	-91.51196354	\\
10	&	145	&	385	&	41.68153099	&	-91.57331503	&	41.64128623	&	-91.56508335	\\
11	&	160	&	340	&	41.67950253	&	-91.57390402	&	41.71443676	&	-91.58289955	\\
12	&	170	&	350	&	41.63511630	&	-91.51468220	&	41.67312935	&	-91.57551051	\\
13	&	190	&	430	&	41.65263384	&	-91.58594751	&	41.65001141	&	-91.46857959	\\
14	&	220	&	400	&	41.65345535	&	-91.52692116	&	41.66699952	&	-91.48110701	\\
15	&	250	&	370	&	41.64186230	&	-91.56777907	&	41.63917026	&	-91.51290389	\\
16	&	255	&	495	&	41.65787087	&	-91.46407211	&	41.64885944	&	-91.55320966	\\
17	&	300	&	420	&	41.70314675	&	-91.60940027	&	41.63583000	&	-91.51710000	\\
18	&	310	&	390	&	41.69834515	&	-91.58892941	&	41.66309000	&	-91.57927000	\\
19	&	330	&	450	&	41.69986134	&	-91.56992240	&	41.65371000	&	-91.49574000	\\
20	&	360	&	580	&	41.65778645	&	-91.56992240	&	41.65529000	&	-91.53254000	\\
21	&	375	&	535	&	41.69835544	&	-91.58876611	&	41.65430000	&	-91.54275000	\\
22	&	390	&	540	&	41.70713817	&	-91.58440115	&	41.66103000	&	-91.54609000	\\
23	&	420	&	620	&	41.61927800	&	-91.53541100	&	41.65157319	&	-91.48818436	\\
24	&	480	&	630	&	41.64975200	&	-91.51395990	&	41.66605613	&	-91.51510378	\\
25	&	525	&	660	&	41.66704200	&	-91.53342870	&	41.68642773	&	-91.51032070	\\
26	&	555	&	705	&	41.64885944	&	-91.55320966	&	41.63202532	&	-91.50501068	\\
27	&	570	&	720	&	41.64199910	&	-91.52728740	&	41.69894765	&	-91.50420625	\\
28	&	590	&	750	&	41.65911430	&	-91.54442830	&	41.64894115	&	-91.58800303	\\
 \hline
\end{tabular}
\end{table}
}
\subsection{Time-Dependent Travel Time}

\noindent For travel beginning at location $i$ at time $t$ and ending at location $j$ within the same time period $w$ while staying within the same region $r$, the speed profile of Table~\ref{table:speedprofile} applies directly and the travel time between $i$ and $j$ is $\tau(i,j,t) = d_{ij} / v_{rw}$, where $d_{ij}$ is the distance (km) between location $i$ and location $j$. For travel between locations $i$ and $j$ in two different regions $r_i$ and $r_j$, and/or across time periods, we describe how we approximate the speed profile in \ref{app:datasets}. Finally, we note it is possible to incorporate other time-dependent or region-dependent factors within the $\tau(i,j,t)$ term, such as expected delays at a pickup location $i$ due to a request not being ready at its quoted ready time or expected delays at delivery location $j$ due to congestion or difficulty locating an address in dense residential areas. 

For travel beginning at location $i$ at time $t$ in time period $w$ and ending at location $j$ such that the travel crosses multiple regions, $r_i \neq r_j$, we approximate the speed profile by weighting the speed profiles of each region intersected by this travel by the size of the respective region. Formally, the travel speed between region $r_i$ and $r_j$ within time period $w$ is approximated by
\begin{eqnarray} \label{eqn:tdtravel}
v_{r_i : r_j w} = \sum_{k \in A(r_i, r_j)} \frac{\alpha_{r_k}}{\sum_{k \in A(r_i, r_j)} \alpha_{r_k}} v_{r_{k} w }  & \forall w,
\end{eqnarray}
where $A(r_i, r_j)$ is the set of regions intersected by traveling between a location $i$ in region $r_i$ to a location $j$ in region $r_j$ and $\alpha_{r}$ is the area of region $r$. For our computational experiments, we apply Equation~(\ref{eqn:tdtravel}) to compute the speed profiles for all $8 \choose 2$ = 28 possible combinations of inter-region travel. 

Let $W$ be the number of time periods in a speed profile. Let the time interval $\left[ l_w, u_w \right]$ define the time period $w$. 
For travel beginning at location $i$ at time $t$ in time period $w$ and ending at location $j$ in time period $w^{\prime} \neq w$, we compute travel time as
\begin{equation} \label{eqn:multtime}
\tau(i,j,t) = \sum_{k=w}^{W} \textnormal{min}\left\{ u_k - \textnormal{max} \left\{ l_k, t \right\}, \frac{d_{ij} - \sum_{s=w}^{k-1} (u_s - \textnormal{max} \left\{ l_s, t \right\}) v_{rs}}{v_{rk}} \right\}.
\end{equation}
If travel from location $i$ to location $j$ intersects multiple regions, then $v_{rk}$ is replaced by $v_{r_i :r_j w}$ in Equation~(\ref{eqn:multtime}).


\section{Additional Figures from Computational Results}
\label{app:morefigures}
\begin{figure}[htb]
\caption{Percent reduction in average cost per request by DRACE versus myopic ALNS.}
\label{figura:diffavgtotalcost2}
\begin{center}
\includegraphics[height=6.25cm,width=6.25cm,angle=0]{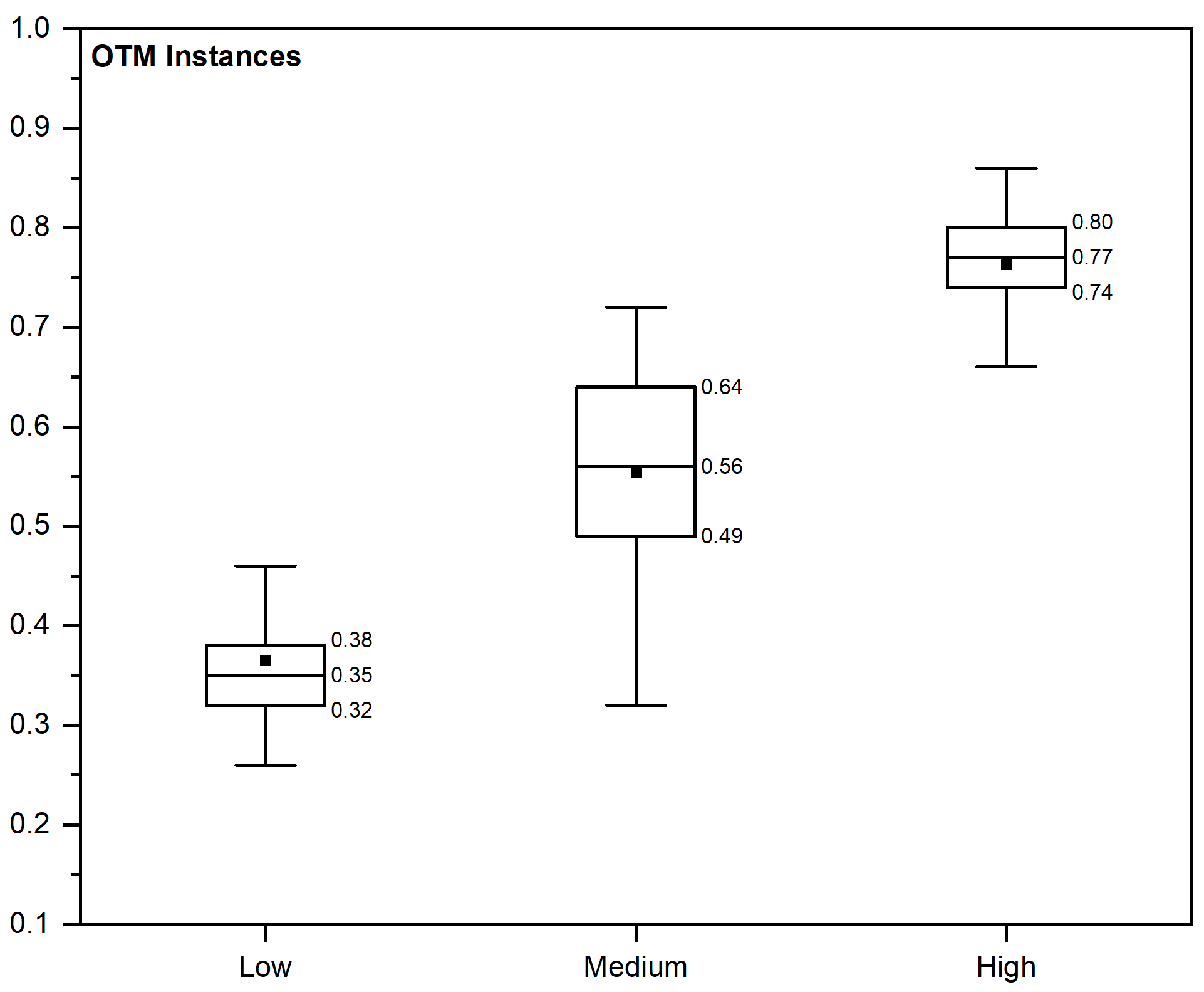}
\includegraphics[height=6.25cm,width=6.25cm,angle=0]{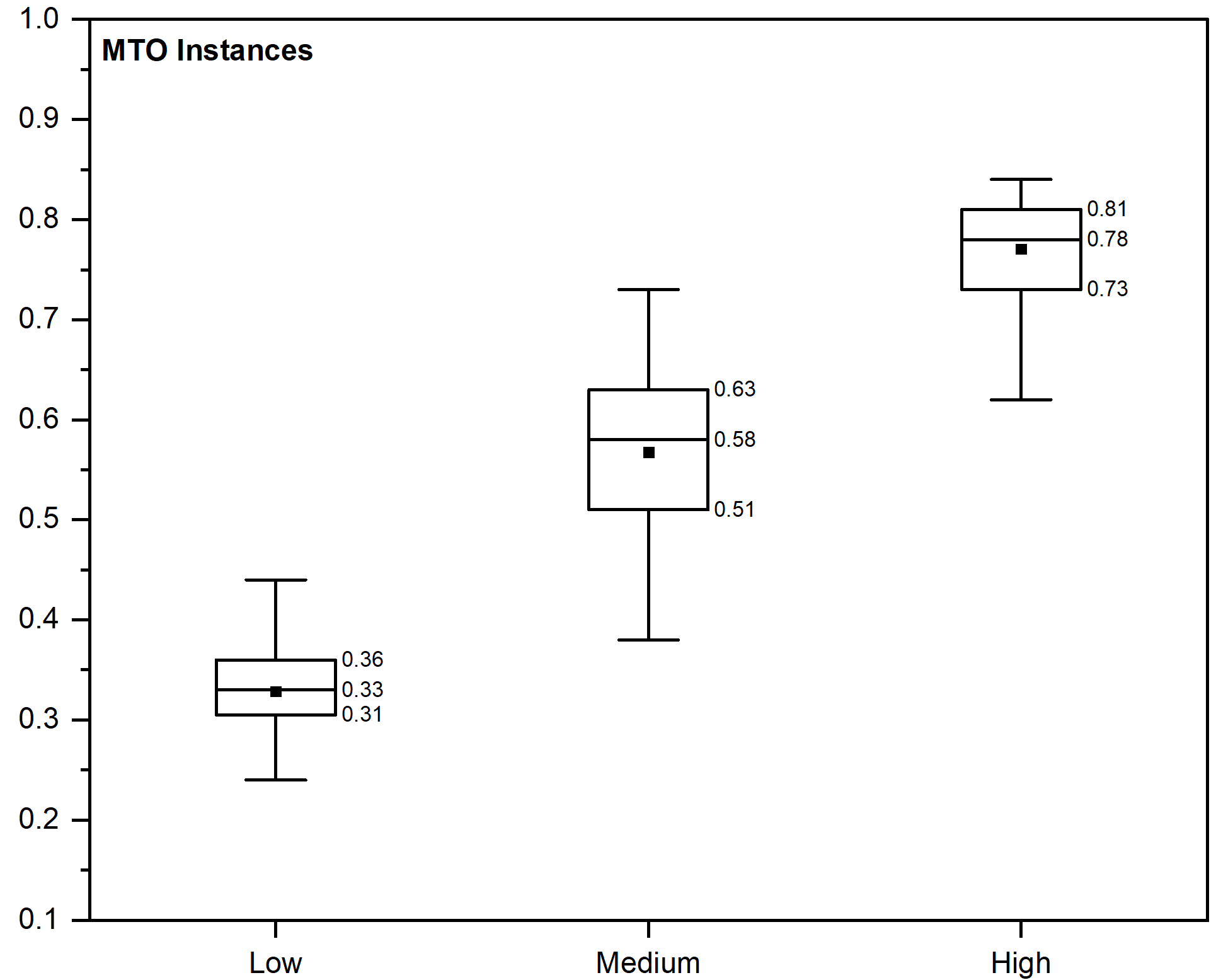}
\end{center}
\end{figure}

\begin{figure}[htb]
\caption{Percent reduction in average lateness per request by DRACE versus myopic ALNS. }   
\label{figura:diffAvgDelay12}
\begin{center}
\includegraphics[height=6.25cm,width=6.25cm,angle=0]{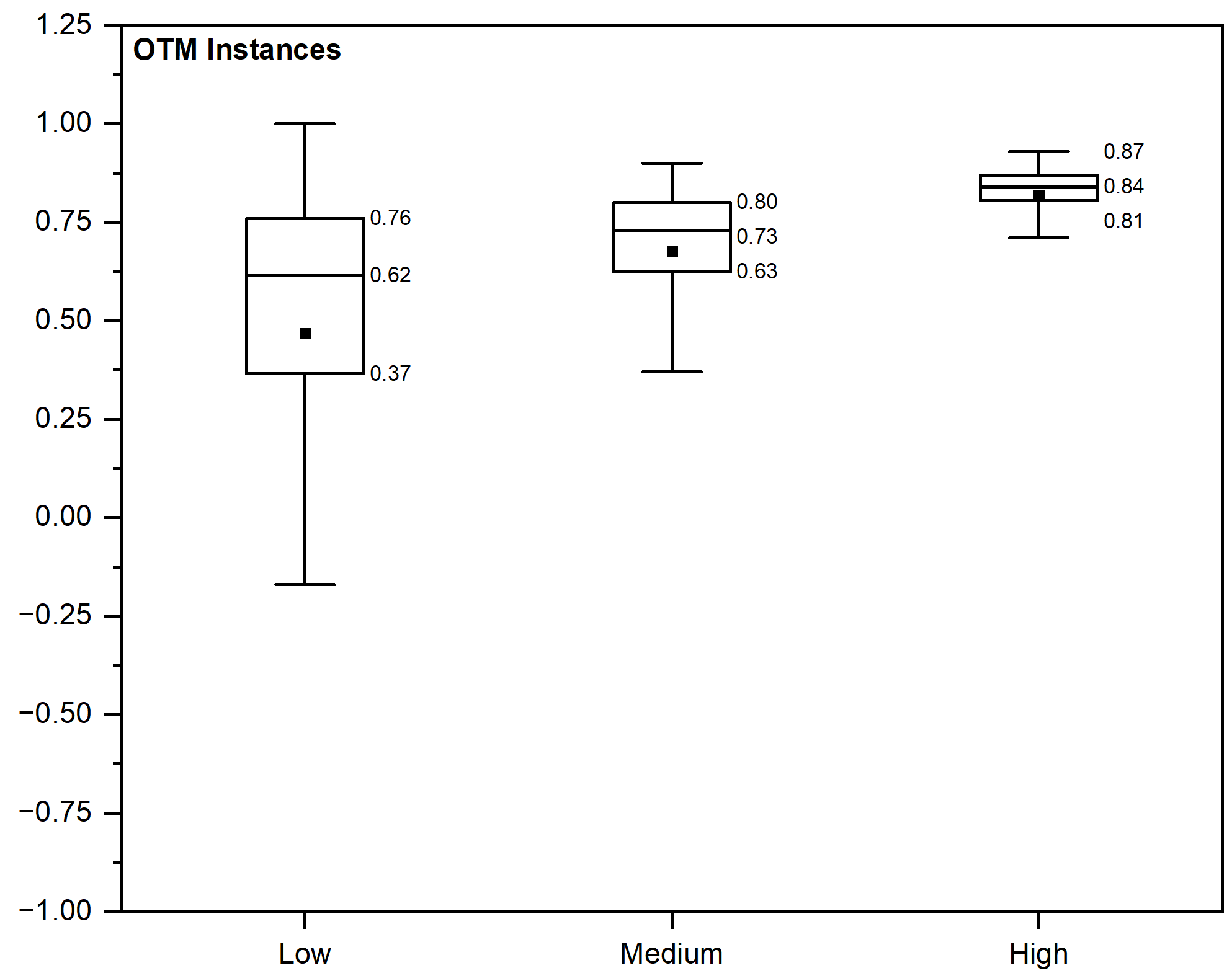}
\includegraphics[height=6.25cm,width=6.25cm,angle=0]{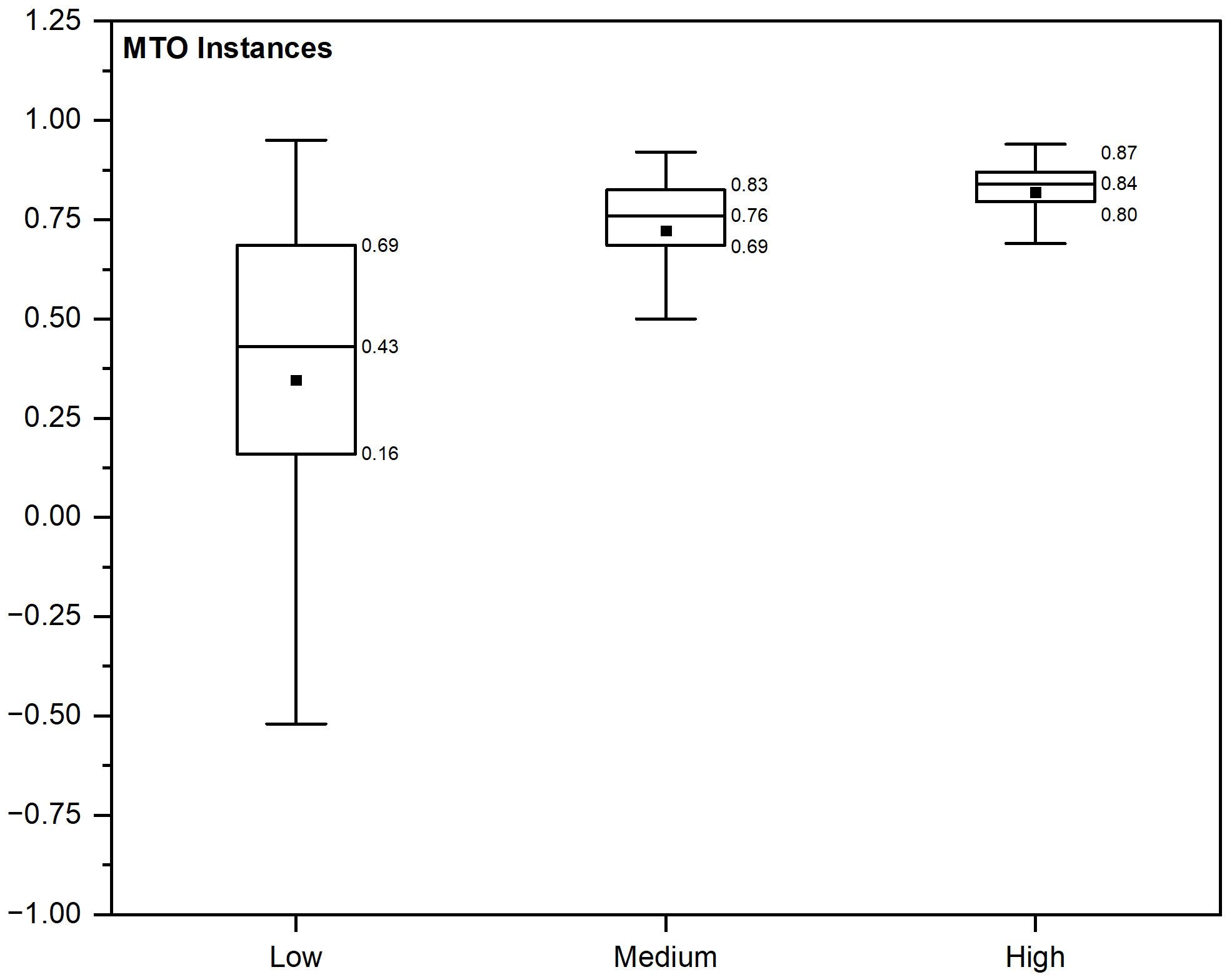}
\end{center}
\end{figure}

\begin{figure}[!h]
\caption{Percent reduction in average cost per request with strategic waiting versus no strategic waiting. }
\label{figura:diffwt_cfa_mto2}


\begin{center}
\includegraphics[height=6.25cm,width=6.25cm,angle=0]{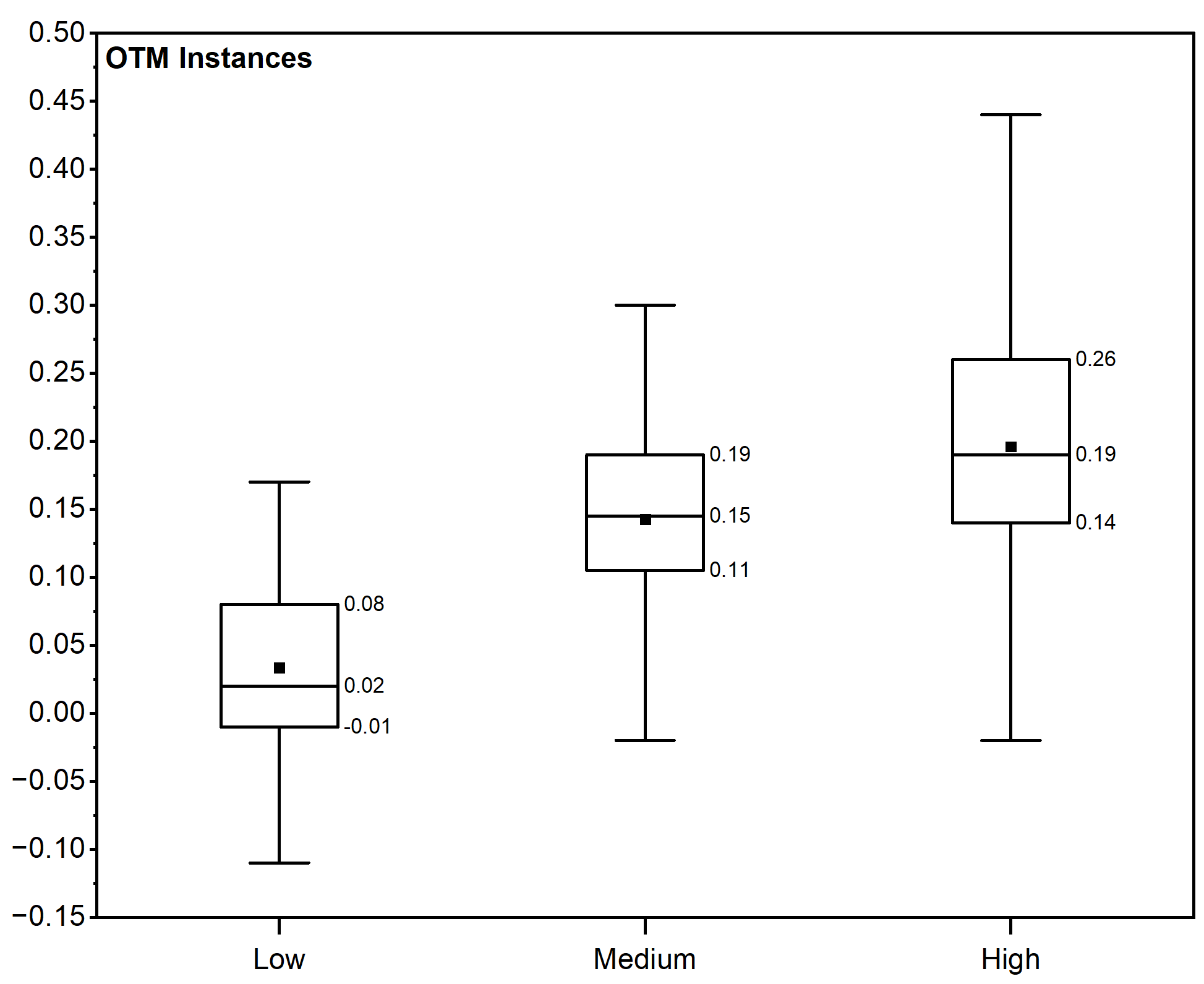}
\includegraphics[height=6.25cm,width=6.25cm,angle=0]{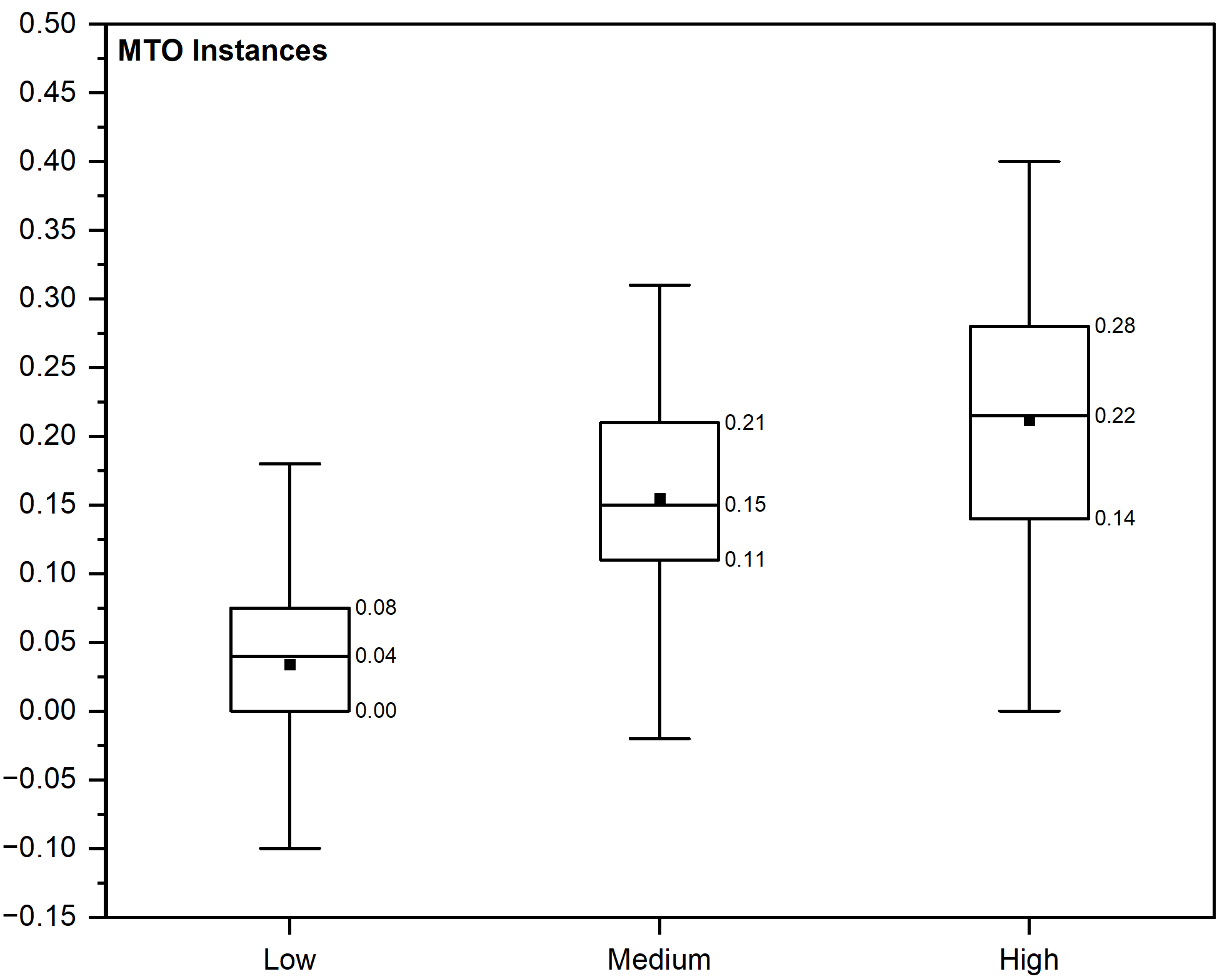}
\end{center}
\end{figure}

\section{Extended Computational Results}
\label{app:moreresults}

\noindent This section comprehensively summarizes all our computational experiments. For each value in the following tables, we calculate the statistics on the output of 100 instances of each type/demand setting. The measure \textit{Total cost} refers to the sum of routing costs and lateness charges for the requests served in a daily instance. \textit{Routing cost} is the sum of travel-related costs and crowdshipper fees for the requests served in a daily instance. \textit{Lateness charge} is the total penalized cost for late deliveries. \textit{Delayed requests} refers to the number of requests served after their deadline. \textit{Total delay} is the total amount of minutes of delay for deliveries over a daily instance. \textit{Crowd service} refers to the number of requests served by the crowdsourced fleet during a day. \textit{Dedicated service} is the number of requests served by the dedicated fleet during a day.

\begin{table}[b]
    \centering
    \begin{tabular}
    {l l r r r r r}
    \hline
        Low NM & KPI & Average & Median & St. Dev.  & Min & Max \\ \hline
          & Total cost & 4867 & 4166 & 375 & 3166 & 5496 \\ 
        ~ & Routing cost & 4722 & 4058 & 271 & 3165 & 4670 \\ 
           & Lateness charge & 145 & 52 & 263 & 1 & 1859 \\ 
        Myopic & Delayed requests & 3 & 2 & 3 & 0 & 25 \\ 
        ~ & Total delay & 8 & 7 & 0 & 0 & 54 \\ 
        ~ & Crowd service & 103 & 89 & 10 & 69 & 119 \\ 
        ~ & Dedicated service & 154 & 133 & 62 & 107 & 152 \\ \hline
         & Total cost & 2819 & 2818 & 223 & 2309 & 3710 \\ 
        ~ & Routing cost & 2686 & 2706 & 170 & 2304 & 3086 \\ 
         & Lateness charge & 133 & 110 & 108 & 5 & 715 \\ 
        DRACE & Delayed requests  & 6 & 6 & 3 & 1 & 15 \\ 
        ~ & Total delay & 4 & 4 & 0 & 0 & 10 \\
        ~ & Crowd service & 206 & 205 & 14 & 170 & 242 \\ 
        ~ & Dedicated service & 15 & 15 & 5 & 3 & 28 \\ \hline
    \end{tabular}
\end{table}

\begin{table}[b]
    \centering
    \begin{tabular}
    {l l r r r r r}
    \hline
        Medium NM & KPI & Average & Median & St. Dev.  & Min & Max \\ \hline
          & Total cost & 11803 & 11086 & 3476 & 5563 & 20575 \\ 
        ~ & Routing cost & 4771 & 4775 & 239 & 4289 & 5391 \\ 
        ~ & Lateness charge & 7032 & 6529 & 3331 & 1217 & 15470 \\ 
        Myopic & Delayed requests & 45 & 45 & 9 & 24 & 67 \\ 
        ~ & Total delay & 30 & 29 & 0 & 0 & 56 \\ 
        ~ & Crowd service & 111 & 113 & 8 & 84 & 128 \\ 
        ~ & Dedicated service & 183 & 185 & 13 & 146 & 208 \\ \hline
         & Total cost & 4461 & 4480 & 410 & 3578 & 5883 \\ 
        ~ & Routing cost & 3950 & 3957 & 271 & 3473 & 4503 \\ 
        ~ & Lateness charge & 511 & 485 & 241 & 75 & 1380 \\ 
        DRACE & Delayed requests & 15 & 15 & 5 & 6 & 31 \\ 
        ~ & Total delay & 7 & 6 & 0 & 0 & 18 \\ 
        ~ & Crowd service & 238 & 238 & 12 & 207 & 267 \\ 
        ~ & Dedicated service & 57 & 57 & 11 & 33 & 78 \\ \hline
    \end{tabular}
\end{table}

\begin{table}[]
    \centering
    \begin{tabular}
    {l l r r r r r}
    \hline
        High NM & KPI & Average & Median & St. Dev.  & Min & Max \\ \hline
         & Total cost & 27804 & 26345 & 7329 & 13001 & 48971 \\ 
        ~ & Routing cost & 5728 & 5722 & 7078 & 5006 & 6523 \\ 
        ~ & Lateness charge & 22076 & 20789 & 7078 & 7723 & 42449 \\ 
        Myopic & Delayed requests & 84 & 81 & 16 & 55 & 137 \\ 
        ~ & Total delay  & 51 & 51 & 0 & 0 & 70 \\ 
        ~ & Crowd service & 141 & 141 & 9 & 122 & 170 \\ 
        ~ & Dedicated service & 230 & 229 & 14 & 199 & 267 \\ \hline
         & Total cost & 5943 & 5828 & 739 & 4601 & 7819 \\ 
        ~ & Routing cost & 4783 & 4801 & 219 & 4236 & 5251 \\
        ~ & Lateness charge & 1160 & 1123 & 617 & 300 & 2975 \\ 
        DRACE & Delayed requests & 26 & 26 & 8 & 10 & 49 \\ 
        ~ & Total delay  & 8 & 8 & 0 & 0 & 17 \\ 
        ~ & Crowd service & 285 & 284 & 15 & 241 & 334 \\ 
        ~ & Dedicated service & 87 & 87 & 11 & 57 & 120 \\ \hline
    \end{tabular}
\end{table}

\begin{table}[]
    \centering
    \begin{tabular}
    {l l r r r r r}
    \hline
        Low OTM   & KPI & Average & Median & St. Dev.  & Min & Max \\ \hline
         & Total cost & 4558 & 4075 & 1530 & 3245 & 5064 \\ 
        ~ & Routing cost & 4084 & 4075 & 309 & 3245 & 5064 \\ 
        ~ & Lateness charge & 474 & 70 & 1347 & 0 & 8428 \\ 
        Myopic & Delayed requests & 6 & 2 & 10 & 0 & 51 \\ 
        ~ & Total delay  & 9 & 8 & 0 & 0 & 33 \\ 
        ~ & Crowd service & 89 & 132 & 9 & 67 & 118 \\ 
        ~ & Dedicated service & 132 & 133 & 9 & 109 & 153 \\ \hline
         & Total cost & 2740 & 2759 & 224 & 2230 & 3442 \\ 
        ~ & Routing cost & 2644 & 2661 & 175 & 2205 & 3122 \\ 
        ~ & Lateness charge & 96 & 70 & 89 & 0 & 585 \\ 
        DRACE & Delayed requests & 5 & 5 & 3 & 0 & 17 \\ 
        ~ & Total delay  & 4 & 3 & 0 & 0 & 17 \\ 
        ~ & Crowd service & 206 & 205 & 14 & 169 & 242 \\ 
        ~ & Dedicated service & 16 & 15 & 5 & 6 & 30 \\ \hline
    \end{tabular}
\end{table}

\begin{table}[]
    \centering
    \begin{tabular}
    {l l r r r r r}
    \hline
        Medium OTM  & KPI & Average & Median & St. Dev.  & Min & Max \\ \hline
         & Total cost & 10633 & 10169 & 2810 & 5922 & 17948 \\ 
        ~ & Routing cost & 4808 & 4772 & 260 & 4203 & 5524 \\ 
        ~ & Lateness charge & 5825 & 5261 & 2630 & 1503 & 12561 \\ 
        Myopic & Delayed requests & 41 & 40 & 9 & 26 & 67 \\ 
        ~ & Total delay  & 27 & 27 & 0 & 0 & 43 \\ 
        ~ & Crowd service & 114 & 115 & 8 & 87 & 132 \\ 
        ~ & Dedicated service & 180 & 180 & 13 & 147 & 205 \\ \hline
         & Total cost & 4496 & 4453 & 442 & 3646 & 6310 \\ 
        ~ & Routing cost & 4007 & 3998 & 240 & 3517 & 4570 \\ 
        ~ & Lateness charge & 489 & 423 & 289 & 85 & 1860 \\ 
        DRACE & Delayed requests & 13 & 13 & 5 & 3 & 26 \\ 
        ~ & Total delay  & 8 & 7 & 0 & 0 & 19 \\ 
        ~ & Crowd service & 236 & 236 & 12 & 211 & 266 \\ 
        ~ & Dedicated service & 58 & 57 & 11 & 33 & 83 \\ \hline
    \end{tabular}
\end{table}

\begin{table}[]
    \centering
    \begin{tabular}
    {l l r r r r r}
    \hline
        High OTM & KPI & Average & Median & St. Dev.  & Min & Max \\ \hline
          & Total cost & 25870 & 24486 & 6653 & 12217 & 45908 \\ 
        ~ & Routing cost & 5733 & 5716 & 321 & 4528 & 6595 \\ 
        ~ & Lateness charge & 20137 & 18754 & 6442 & 6937 & 39465 \\ 
         Myopic & Delayed requests & 80 & 77 & 15 & 51 & 125 \\ 
        ~ & Total delay  & 49 & 50 & 0 & 0 & 65 \\ 
        ~ & Crowd service & 143 & 142 & 9 & 125 & 172 \\ 
        ~ & Dedicated service & 229 & 228 & 14 & 198 & 265 \\ \hline
         & Total cost & 5834 & 5711 & 835 & 4469 & 10831 \\ 
        ~ & Routing cost & 4772 & 4763 & 231 & 4251 & 5367 \\ 
        ~ & Lateness charge & 1062 & 915 & 708 & 195 & 5585 \\ 
        DRACE & Delayed requests & 24 & 24 & 8 & 9 & 58 \\ 
        ~ & Total delay  & 8 & 8 & 0 & 0 & 19 \\ 
        ~ & Crowd service & 285 & 284 & 13 & 259 & 323 \\ 
        ~ & Dedicated service & 87 & 85 & 12 & 62 & 119 \\ \hline
    \end{tabular}
\end{table}

\begin{table}[]
    \centering
    \begin{tabular}
    {l l r r r r r}
    \hline
        Low MTO  & KPI & Average & Median & St. Dev.  & Min & Max \\ \hline 
        ~ & Total cost & 4070 & 4109 & 299 & 3283 & 5345 \\ 
        ~ & Routing cost & 3973 & 3972 & 251 & 3281 & 4540 \\ 
        ~ & Lateness charge & 97 & 50 & 134 & 2 & 854 \\ 
        Myopic & Delayed requests & 2 & 2 & 2 & 0 & 9 \\ 
        ~ & Total delay  & 8 & 6 & 0 & 0 & 47 \\ 
        ~ & Crowd service & 87 & 87 & 10 & 63 & 116 \\ 
        ~ & Dedicated service & 135 & 135 & 8 & 111 & 154 \\ \hline
        ~ & Total cost & 2731 & 2750 & 237 & 2170 & 3360 \\ 
        ~ & Routing cost & 2611 & 2620 & 181 & 2150 & 3034 \\ 
        ~ & Lateness charge & 120 & 100 & 99 & 10 & 485 \\ 
        DRACE & Delayed requests & 6 & 6 & 3 & 1 & 14 \\ 
        ~ & Total delay  & 4 & 3 & 0 & 0 & 12 \\ 
        ~ & Crowd service & 206 & 205 & 14 & 172 & 240 \\ 
        ~ & Dedicated service & 15 & 14 & 5 & 4 & 30 \\ \hline
    \end{tabular}
\end{table}

\begin{table}[]
    \centering
    \begin{tabular}
    {l l r r r r r}
    \hline
       Medium MTO  & KPI & Average & Median & St. Dev.  & Min & Max \\ \hline
        ~ & Total cost & 10885 & 10750 & 2948 & 5458 & 19588 \\ 
        ~ & Routing cost & 4716 & 4739 & 257 & 4161 & 5244 \\ 
        ~ & Lateness charge & 6169 & 6158 & 2782 & 4161 & 5244 \\ 
        Myopic & Delayed requests & 43 & 43 & 9 & 22 & 67 \\ 
        ~ & Total delay  & 27 & 27 & 0 & 0 & 44 \\ 
        ~ & Crowd service & 112 & 112 & 7 & 95 & 129 \\ 
        ~ & Dedicated service & 182 & 185 & 13 & 147 & 208 \\ \hline
        ~ & Total cost & 4459 & 4461 & 466 & 3578 & 6353 \\ 
        ~ & Routing cost & 3909 & 3933 & 253 & 3468 & 4484 \\ 
        ~ & Lateness charge & 550 & 538 & 289 & 75 & 2150 \\ 
        DRACE & Delayed requests & 16 & 16 & 5 & 6 & 31 \\ 
        ~ & Total delay  & 7 & 7 & 0 & 0 & 18 \\ 
        ~ & Crowd service & 238 & 238 & 12 & 212 & 265 \\ 
        ~ & Dedicated service & 56 & 55 & 11 & 33 & 80 \\ \hline
    \end{tabular}
\end{table}

\begin{table}[]
    \centering
    \begin{tabular}
    {l l r r r r r}
    \hline
        High MTO & KPI & Average & Median & St. Dev.  & Min & Max \\ \hline
        ~ & Total cost & 26724 & 25561 & 7116 & 13242 & 47522 \\ 
        ~ & Routing cost & 5823 & 5646 & 1102 & 4956 & 14959 \\ 
        ~ & Lateness charge & 20901 & 19891 & 6924 & 8287 & 40990 \\ 
        Myopic & Delayed requests & 83 & 80 & 16 & 52 & 137 \\ 
        ~ & Total delay  & 49 & 50 & 0 & 0 & 64 \\ 
        ~ & Crowd service & 141 & 140 & 9 & 116 & 174 \\ 
        ~ & Dedicated service & 231 & 230 & 13 & 198 & 264 \\ \hline
        ~ & Total cost & 5822 & 5571 & 798 & 4569 & 9495 \\ 
        ~ & Routing cost & 4665 & 4635 & 238 & 3955 & 5347 \\ 
        ~ & Lateness charge & 1157 & 1025 & 658 & 305 & 4265 \\ 
        DRACE & Delayed requests & 26 & 26 & 8 & 10 & 54 \\ 
        ~ & Total delay  & 8 & 8 & 0 & 0 & 16 \\ 
        ~ & Crowd service & 288 & 285 & 15 & 258 & 330 \\ 
        ~ & Dedicated service & 84 & 83 & 12 & 57 & 119 \\ \hline
    \end{tabular}
\end{table}

\begin{table}[]
    \centering
    \begin{tabular}
    {l l r r r r r}
    \hline
        Low UO & KPI & Average & Median & St. Dev.  & Min & Max \\\hline 
        ~ & Total cost & 14670 & 13541 & 5296 & 6590 & 28656 \\ 
        ~ & Routing cost & 2414 & 2383 & 203 & 2019 & 2874 \\ 
        ~ & Lateness charge & 12256 & 11099 & 5122 & 4481 & 25783 \\ 
        Myopic & Delayed requests & 58 & 55 & 14 & 35 & 91 \\ 
        ~ & Total delay  & 41 & 40 & 8 & 24 & 58 \\ 
        ~ & Crowd service & 69 & 68 & 10 & 52 & 94 \\ 
        ~ & Dedicated service & 109 & 109 & 6 & 96 & 121 \\ \hline
        ~ & Total cost & 2895 & 2783 & 603 & 2109 & 4855 \\ 
        ~ & Routing cost & 1981 & 1979 & 116 & 1679 & 2335 \\ 
        ~ & Lateness charge & 914 & 793 & 549 & 230 & 2880 \\ 
        DRACE & Delayed requests & 16 & 16 & 5 & 7 & 33 \\ 
        ~ & Total delay  & 11 & 10 & 5 & 4 & 36 \\ 
        ~ & Crowd service & 161 & 160 & 12 & 137 & 189 \\ 
        ~ & Dedicated service & 19 & 18 & 4 & 10 & 32 \\ \hline
    \end{tabular}
\end{table}

\begin{table}[]
    \centering
    \begin{tabular}
    {l l r r r r r}
    \hline
        Medium UO & KPI & Average & Median & St. Dev.  & Min & Max \\ \hline
        ~ & Total cost & 39423 & 38173 & 9540 & 19465 & 83466 \\ 
        ~ & Routing cost & 3733 & 3318 & 2503 & 2555 & 21346 \\ 
        ~ & Lateness charge & 35690 & 34510 & 9905 & 12214 & 79236 \\ 
        Myopic & Delayed requests & 120 & 118 & 21 & 68 & 199 \\ 
        ~ & Total delay  & 59 & 59 & 8 & 25 & 80 \\ 
        ~ & Crowd service & 114 & 114 & 14 & 77 & 157 \\ 
        ~ & Dedicated service & 127 & 127 & 6 & 112 & 142 \\ \hline
        ~ & Total cost & 6390 & 6180 & 1555 & 3395 & 11221 \\ 
        ~ & Routing cost & 2510 & 2507 & 140 & 2159 & 2823 \\ 
        ~ & Lateness charge & 3881 & 3655 & 1499 & 985 & 8510 \\ 
        DRACE & Delayed requests & 41 & 38 & 10 & 20 & 69 \\ 
        ~ & Total delay  & 19 & 19 & 9 & 9 & 30 \\ 
        ~ & Crowd service & 211 & 210 & 16 & 161 & 253 \\ 
        ~ & Dedicated service & 32 & 32 & 6 & 20 & 52 \\ \hline
    \end{tabular}
\end{table}

\begin{table}[]
    \centering
    \begin{tabular}
    {l l r r r r r}
    \hline
        High UO  & KPI & Average & Median & St. Dev.  & Min & Max \\ \hline
        ~ & Total cost & 66743 & 66188 & 9261 & 45013 & 93327 \\ 
        ~ & Routing cost & 4710 & 4041 & 6899 & 3356 & 72976 \\ 
        ~ & Lateness charge & 62033 & 62036 & 10999 & 355 & 88825 \\ 
        Myopic & Delayed requests & 180 & 181 & 15 & 135 & 216 \\ 
        ~ & Total delay  & 69 & 68 & 5 & 59 & 88 \\ 
        ~ & Crowd service & 152 & 152 & 11 & 127 & 181 \\ 
        ~ & Dedicated service & 137 & 138 & 4 & 120 & 146 \\ \hline
        ~ & Total cost & 9838 & 9672 & 2263 & 4397 & 16284 \\ 
        ~ & Routing cost & 2896 & 2906 & 155 & 2416 & 3308 \\ 
        ~ & Lateness charge & 6942 & 6745 & 2218 & 1580 & 13390 \\ 
        DRACE & Delayed requests & 63 & 62 & 12 & 38 & 101 \\ 
        ~ & Total delay  & 22 & 22 & 8 & 8 & 34 \\ 
        ~ & Crowd service & 256 & 256 & 12 & 227 & 285 \\ 
        ~ & Dedicated service & 45 & 45 & 7 & 30 & 68 \\ \hline
    \end{tabular}
\end{table}

\end{document}